\documentclass[11pt]{article}
\usepackage{latexsym}
\usepackage{amsmath,amssymb}
\usepackage{amsthm}
\usepackage{stmaryrd}
\usepackage{float}
\usepackage{caption}
\title{Thompson-like characterization of solubility for products of finite groups\thanks{ Research supported by Proyectos PROMETEO/2017/057 from the Generalitat Valenciana (Valencian Community, Spain), and PGC2018-096872-B-I00 from the Ministerio de Ciencia, Innovaci\'on y Universidades, Spain, and FEDER, European
Union; and second author also by Project VIP-008 of Yaroslavl P. Demidov State University.}}
\author{P. Hauck, L.~S. Kazarin,  A. Mart\'{\i}nez-Pastor, M.~D. P\'{e}rez-Ramos\\
\\
\\
{\it Dedicated to the memory of M.~Pilar~G\'{a}llego (1958-2019)}}

\newcommand{\N}{\mathcal{N}}

\newcommand{\So}{\mathcal{S}}

\newtheorem{lem}{Lemma}
\newtheorem{pro}{Proposition}
\newtheorem{te}{Theorem}
\newtheorem*{teo}{Theorem}
\newtheorem{cor}{Corollary}
\theoremstyle{definition}

\newtheorem{de}{Definition}
\newtheorem*{notation}{Notation}
\newtheorem{primes}{Notation}
\newtheorem*{Strat}{Strategies}
\theoremstyle{remark}
\newtheorem*{rem}{Remark}

\begin{document}
\date{}
\maketitle

\begin{abstract}
A remarkable result of Thompson states that a finite group is soluble if and only if its two-generated subgroups are soluble. This result has been generalized in numerous ways, and it is in the core of a wide area of research in the theory of groups, aiming for global properties of groups from local properties of two-generated (or more generally, $n$-generated) subgroups. We contribute an extension of Thompson's theorem from the perspective of factorized groups. More precisely, we study finite groups $G = AB$ with subgroups $A,\ B$ such that $\langle a, b\rangle$ is soluble for all $a \in A$ and $b \in B$. In this case, the group $G$ is said
to be an {\it $\cal S$-connected product} of the subgroups $A$ and $B$ for the class $\cal S$ of all finite soluble groups. Our main theorem states that $G = AB$ is $\cal S$-connected
if and only if $[A,B]$ is soluble. In the course of the proof we derive a result
of own interest about independent primes regarding the soluble graph of
almost simple groups.\medskip

\medskip

\end{abstract}

\section{Introduction}

Our work arises from the confluence of two major areas of study in group theory. On the one hand, what we might call local-global theory, and on the other, the theory of products of groups.
\smallskip

Regarding the first one, the interest lies in the influence on the global structure of a group of local properties on its elements, either by satisfying explicit relations or formulas or by their generation properties. Classical Burnside problems might be traced to the origin of this theory.  We paraphrase here F.~Grunewald, B.~Kunyavski\u{\i} and E.~Plotkin in \cite{GKP}, which provides a valuable reference on the topic. Widely speaking, one can say that  classical Burnside problems ask to what extent finiteness  of cyclic subgroups (i.e. generated by one element) determines finiteness of arbitrary finitely generated subgroups of a group.  We are interested in the influence of two-generated subgroups on the structure of finite groups. First results in
this direction by M.~Zorn \cite{Zo} and R.~Baer \cite{B} show that nilpotency and
supersolubility, respectively, of a finite group is determined by the same corresponding property of its two-generated subgroups.  Undoubtedly one of the most influential results is the one of J.~Thompson   regarding solubility.
\begin{teo}\textup{(Thompson,  \cite{T})} A finite group $G$ is soluble if and only if every two-generated subgroup of $G$ is soluble.
\end{teo}
This result has been generalized and sharpened in various
ways. In addition to the above-mentioned reference,  we cite for instance \cite{DGHP, GGKP1, GGKP2, Gu, GL, GKPS, MSW}, some of whose results have been applied to prove the results in this paper. As a typical example we mention the following theorem of R.~Guralnick, K.~Kunyavski\u{\i},
E.~Plotkin and A.~Shalev.
\begin{teo}\textup{(Guralnick, Kunyavski\u{\i},
Plotkin, Shalev,  \cite{GKPS})} Let $G$ be a finite group, let $G_{\cal S}$ denote the soluble radical of $G$ (i.e. the largest soluble normal subgroup of $G$) and let $x\in G$. Then $x\in G_{\cal S}$ if and only if the subgroup $\langle x,y\rangle$ is soluble for all $y\in G$. 
\end{teo}
\smallskip

In the theory of products of groups, the aim is to seek for information about the structure of a factorized group from the subgroups in the factorization (and vice versa). The well-known theorem by Kegel and Wielandt, about the solubility of a finite group which is the product of nilpotent subgroups, is probably one of the most remarkable results in the area. It is also known that the product of two finite supersoluble subgroups is not necessarily supersoluble, even if the factors are normal subgroups.
This fact has motivated the search for conditions to obtain positive results and, at the time, has been the source of a vast line of research on factorized finite groups whose factors are linked by some particular property. Originally M.~Asaad and A.~Shaalan in \cite{AS} introduced totally permutable products of subgroups, which can be seen as extension of central products. A group $G=AB$ is a central product of subgroups $A$ and $B$ if $ab=ba$, for all $a\in A$ and $b\in B$; equivalently, $\langle a,b\rangle$ is abelian,  for all $a\in A$ and $b\in B$. The subgroups $A$ and $B$ are said to be totally permutable if  every subgroup of $A$ permutes with every subgroup of $B$; equivalently, $\langle a\rangle\langle b\rangle= \langle b\rangle\langle a\rangle$, for all $a\in A$ and $b\in B$. R.~Maier notes in \cite{M} that for such subgroups, $\langle a,b\rangle=\langle a\rangle\langle b\rangle=\langle b\rangle\langle a\rangle$ is supersoluble, for all $a\in A$ and $b\in B$, which led to the following connection property:
\begin{de}\textup{(Carocca, \cite{C96})} Let $\cal L$ be a non-empty class of groups. Subgroups $A$ and $B$ of a group $G$ are $\cal L$-connected if $\langle a,b\rangle\in \cal L$ for all $a\in A$ and $b\in B$.
\end{de}

For the special case when $G=AB=A=B$ this means of course that $\langle a,b\rangle\in \cal L$ for all $a,b\in G$, and the study of products of $\cal L$-connected subgroups provides a more general setting for local-global questions related to two-generated subgroups. We refer to \cite{Ncon1, Ncon2, Ncon3} for previous studies for the class $\mathcal L=\cal N$ of finite nilpotent groups, and to \cite{GHP1, GHP2, GHP3, GHP4} for $\mathcal L$ being the class of finite metanilpotent groups and other relevant classes of groups. For the class $\mathcal L=\cal S$ of finite soluble groups, A.~Carocca in \cite{C} proved the solubility of a product of $\cal S$-connected soluble subgroups, which provides a first extension of the above-mentioned theorem of Thompson for products of groups (see Corollary~\ref{corCa}).
\smallskip

All groups considered in this paper are assumed to be finite. Unless otherwise specified, we shall adhere to the notation used in \cite{DH} and we refer also to that book for the basic results on classes of groups.
In particular, $G_{\So}$ denotes the soluble radical of a group $G$ as mentioned before.
In addition, if $n$ is a positive integer, then $\pi(n)$ denotes the set of primes dividing $n$;  and $\pi(H) = \pi(|H|)$ for any group $H$.
\smallskip

The main result in this paper is the following:
\smallskip

{\bf Main Theorem.} {\it Let the finite group $G = AB$ be the product of subgroups $A$ and $B$. Then the following statements
are equivalent:
\begin{itemize} \item[(1)] $A, B$ are $\So$-connected.
\item[(2)] For all primes $p \ne q$, all $p$-elements $a \in A$ and all $q$-elements $b \in B$, $\langle a,b \rangle$
is soluble.
\item[(3)] $[A,B] \le G_{\So}$.
\end{itemize} }

The configuration of a minimal counterexample to the Main Theorem is proven to be an almost simple group, i.e. a group $G$ such that $N\unlhd G\le \mbox{\upshape Aut}(N)$ for some non-abelian simple group $N$, where $\mbox{\upshape Aut}(N)$ denotes the automorphism group of $N$. As a major previous result, in Section~2 we prove Theorem~\ref{as}, which deals with an almost simple group which is the product of subgroups satisfying condition~(2) of the Main Theorem. A stronger version of Theorem~\ref{as} is stated in Corollary~\ref{cote1} as a consequence of our  main theorem. A remarkable result deduced from the checking carried out for the proof of Theorem~\ref{as} is stated in Theorem~\ref{fact} and has to do with the existence of {\it independent primes} in almost simple groups.

For a group $G$ the Gr\"unberg-Kegel or prime graph $\Gamma(G)$ of $G$ is well-known.
It consists of the set of vertices $V(\Gamma(G))=\pi(G)$ and the set of edges $(p,q)\in \pi(G)\times \pi(G)$ such that there is an element of order $pq$ in $G$. S.~Abe and N.~Iiyori introduced in \cite{AI} the soluble graph $\Gamma_{sol}(G)$ which has the same set of vertices as $\Gamma(G)$, but two vertices $p,q\in \pi(G)$ are  adjacent if $G$ contains a soluble subgroup of order divisible by $pq$. Iiyori~\cite{I}, B.~Amberg, A.~Carocca and L.~Kazarin \cite{ACK}, and Amberg and Kazarin \cite{AK} take further the study of the soluble graph specially for non-abelian finite simple groups. In the last reference the authors  are concerned with subsets $I\subseteq \pi(G)$ such that no pair of vertices in $I$ is adjacent with respect to $\Gamma(G)$ or $\Gamma_{sol}(G)$, respectively, called {\it independent sets} of the corresponding graph.  This leads to the following main concept for our purposes.

\begin{de} Given a finite group $H$, we call two prime divisors $p$ and $q$ of $|H|$ {\it independent}
(with respect to $H$), if $H
$ contains no soluble subgroup whose order is divisible by $pq$.
\end{de}

For an almost simple group $G=AB$ with non-abelian socle $N$, not contained in either $A$ or $B$, and except for a few exceptions, we state in Theorem~\ref{fact}, and in the subsequent Remark, the existence of independent primes with respect to $N$, one dividing $|A\cap N|$, the other one dividing $|B\cap N|$, derived from the proof of Theorem~\ref{as}. We also prove that apart from some additional exceptions the existence of such independent primes with  respect to $\mbox{\upshape Aut}(N)$ remains true.
\smallskip

It is clear that the Main Theorem extends the above-mentioned theorems of Thompson~\cite{T} and Carocca~\cite{C}. It also implies the theorem
of Guralnick, Kunyavski\u{\i}, Plotkin and Shalev~\cite{GKPS} stated above (with $A = G$,
$B = \langle x\rangle$; note that $\langle x\rangle G_{\cal S}$ is a normal (soluble) subgroup of G by the Main
Theorem). Section~3 is devoted to prove our main result and to state some first consequences. In particular,  Corollary~\ref{corCa} generalizes Carocca's result via the soluble radical in a product of $\cal S$-connected subgroups. In a forthcoming paper \cite{GHKMP}, our theorem is applied to extend  main results known for finite soluble groups in \cite{GHP1, GHP2, GHP3} to the universe of all finite groups.

\section{The case of almost simple groups}

The aim of this section is the proof of the following result:

\begin{te}\label{as} Let $N$ be a non-abelian simple group and $N \le G = AB \le \mbox{\upshape Aut}(N)$ with
subgroups $A$ and $B$ satisfying condition (2) of the Main Theorem, $AN = BN = G$. Then $A = G$ or $B = G$.
\end{te}

In Corollary~\ref{cote1} we will see, as a consequence of the Main Theorem, that the assumption $AN = BN = G$ in the preceding
result is not necessary and that actually $A = 1$ or $B = 1$ holds.\\

Occasionally, the following lemma will be useful (cf. \cite{AFG}, Lemma 1.3.1).

\begin{lem}\label{con} Let $G = AB$ be a group with subgroups $A, B$ satisfying condition (2) of the Main Theorem. If $g,h \in G$, then $G = A^gB^h$ and $A^g, B^h$ also satisfy condition (2) of the Main Theorem.
\end{lem}

{\bf Proof.} Let $g = a_1b_1$, $h = b_2a_2$, $b_1a_2^{-1} = a_3b_3$ with $a_1,a_2,a_3 \in A$,
$b_1, b_2, b_3 \in B$. If $a \in A$ and $b \in B$, then $a^g = (a^{a_1a_3})^{b_3a_2}$ and
$b^h = (b^{b_2b_3^{-1}})^{b_3a_2}$. The assertion follows.\qed \\

Apart from the case of alternating groups our treatment in this section relies on the classification of
the maximal factorizations of almost simple groups by Liebeck, Praeger and Saxl \cite{LPS} (and \cite{HLS}
for the exceptional groups of Lie type).
\smallskip

In order to prove Theorem~\ref{as}, we assume here that $N$ is neither contained in $A$ nor in $B$
(for otherwise $A = G$ or $B = G$ and we are done). Our aim is to show that for
all possible factorizations $G = AB$, with $N\nleq A$ and $N\nleq B$, the subgroups $A$ and $B$ do not satisfy condition (2) of the Main Theorem.\smallskip

Since we have also
to consider possibly non-maximal factorizations $G = AB$, we use the following notation:

\begin{notation} If $G = AB$ is as above, choose maximal subgroups $\widetilde{A} \ge A$ and $\widetilde{B} \ge B$, $N\nleq \widetilde{A}$ and $N\nleq \widetilde{B}$. Then $G =\widetilde{A}\widetilde{B}$ is a factorization of the type considered in \cite{LPS}.\smallskip

Except for this $\,_{\widetilde{}}$ -notation we use the notation of \cite{LPS}, however with $N$ instead of $L$.
 \end{notation}

We need the following simple lemma which follows from elementary order considerations. To formulate it,
we denote for a positive integer $n$ and a prime $p$ by $n_p$ the highest $p$-power dividing $n$, i.e.
$n = n_pm$, $p \not\:\mid m$.

\begin{lem}\label{l1} Let $G = AB = AN = BN$ as in Theorem \ref{as}, $p$ a prime. Then
$|A \cap N|_p \ge \frac{|N|_p}{|\widetilde{B} \cap N|_p|\mbox{\upshape Out}(N)|_p}$ and likewise with $B$ and
$\widetilde{A}$.
\end{lem}

{\bf Proof.} Since $G=\widetilde{B}N=\widetilde{B}A=AN$, we easily have:
\begin{gather*}  \frac{|N|}{|\widetilde{B} \cap N||\mbox{\upshape Out}(N)|}\le \frac{|N|}{|\widetilde{B} \cap N||G/N|}= \frac{|N|}{|\widetilde{B} \cap N||\widetilde{B}/(\widetilde{B}\cap N)|} = \frac{|N|}{|\widetilde{B}|}\\
=\frac{|N||A|}{|G||\widetilde{B}\cap A|} = \frac{|A\cap N|}{|\widetilde{B}\cap A|}\le |A\cap N|,
\end{gather*}  and the result follows.\qed
\medskip

We state the outline of the proof, and explain general arguments and strategies used to carry it out, though usually detailed checking work and easy calculations are omitted.

\begin{Strat}

\begin{description} \item[(S1)] In our treatment of almost simple groups $G$ (especially those of Lie type and the sporadic ones), we will usually determine
two primes, one dividing $|A \cap N|$
and the other dividing $|B \cap N|$ that are independent with respect to the simple group $N$.  The
divisibility properties follow always from Lemma \ref{l1}.\\
We remark that in sections \ref{class} - \ref{spor} independence of primes is always meant with respect to $N$,
even if not explicitely stated.
\item[(S2)] For certain small groups the independency of two
primes can be deduced from the subgroup structure given in
\cite{BHR} or \cite{Atl}. \item[(S3)] The following observation is sometimes helpful: Given a factorization $G = AB$, suppose
by way of contradiction that $A$, $B$ satisfy condition (2) of the Main Theorem. If $r$ and $s$ are distinct primes,
$a \in A$ an $r$-element and $b \in B$ an $s$-element, we may assume by Lemma \ref{con} that $a$ and $b$ generate an
$\{r,s\}$-subgroup of $G$ by replacing $A$ by a suitable conjugate.

\item[(S4)]For infinite families of groups of Lie type, the independency of the specified primes can usually be
proved by referring to results about the prime divisors of certain maximal soluble subgroups. One of
the primes in question is frequently a primitive prime divisor of $p^k - 1$ for a suitable $k$,
where $p$ is the characteristic corresponding to the group. \end{description} \end{Strat}

We recall this relevant definition:

\begin{de} Let $k$ be a positive integer and $p$ a prime. A {\it primitive prime divisor} of $p^k - 1$ is a prime $r$ such that
$r | p^k - 1$ and $r \not\:\mid p^i -1$ for every integer $i$ such that $1 \le i < k$.\end{de}

The following well-known lemma of Zsigmondy \cite{Z} describes when primitive prime divisors exist:

\begin{lem}\label{Zs}
Let $k \ge 2$ an integer and $p$ a prime.
\begin{itemize}
\item[a)] There exists a primitive prime divisor of $p^k - 1$ unless $k = 2$ and $p$ is a Mersenne prime
or $(p,k) = (2,6)$.
\item[b)] If $r$ is a primitive prime divisor of $p^k - 1$, then $r - 1 \equiv 0\, (\mbox{\upshape mod } k)$. In
particular, $r \ge k + 1$.
\end{itemize}
\end{lem}

\subsection{Alternating groups}

\begin{lem}\label{prime} If $n$ is an integer, $n \ge 5$, $n \neq 10$, then there exists a non-Mersenne prime $p$
such that $\frac{n}{2} < p \le n$.
\end{lem}

{\bf Proof.} By a generalization of the Bertrand-Chebyshev theorem, due to Ramanujan \cite{Ra}, for $n \ge 11$ there exist two
primes $p_1, p_2$ with $\frac{n}{2} < p_1 < p_2 \le n$. Clearly, not both of them can be Mersenne primes. Hence
the assertion holds for $n \ge 11$ and $p = 5$ works for $5 \le n \le 9$.\qed

\begin{lem}\label{pclo} Let $H$ be a soluble permutation group of degree $n$. If there exists a non-Mersenne
prime $p > \frac{n}{2}$ dividing $|H|$, then $H$ is $p$-closed.
\end{lem}

{\bf Proof.} $p^2$ does not divide $|H|$ as $p > \frac{n}{2}$. Hence the assertion is equivalent to
$\langle x \rangle$ normal in $H$ for $x \in H$ of order $p$; note that $x$ is a $p$-cycle.

The proof is by induction on $n$, the case $n = 1$ being trivial. Assume that $H$ is intransitive on $\Omega =
\{1,\ldots,n\}$ with orbits $\Delta_1,\ldots,\Delta_s$. Then $H \le H^{\Delta_1} \times \cdots \times H^{\Delta_s}$
where $H^{\Delta_i}$ is the permutation group induced by $H$ on $\Delta_i$. It follows that there exists exactly
one $j$ with $|\Delta_j| \ge p$ and $p$ divides $|H^{\Delta_j}|$. By induction, $H^{\Delta_j}$ is $p$-closed and
so is $H$.

Hence we may assume that $H$ is transitive on $\Omega$. Then $H$ is primitive on $\Omega$ (see 1.2.(a) in \cite{WW}).
If $n = p$ then $H$ is 2-transitive or $p$-closed by a result of Burnside (see \cite{DM}, Corollary 3.5B). If $n > p$,
then $H$ is 2-transitive by a result of Jordan (see \cite{DM}, Theorem 7.4A). So it remains to consider the
2-transitive case.

According to Huppert's classification of soluble 2-transitive permutation groups \cite{Hup},
$n = r^a$ for a prime $r$ and $H \le AGL_a(r)$ where either $H$ is a group of semilinear mappings over $\mathbb{F}_{r^a}$
with $|H| = r^a(r^a -1)b$, $b|a$ or $n \in \{3^2,3^4,5^2,7^2,11^2,23^2\}$.

It is easily checked that in none of the exceptional cases there exists a prime divisor of $|AGL_a(r)|$ which
is larger than $\frac{r^a}{2}$.

In the generic case, $|H| = r^a(r^a -1)b$, $b|a$. Since $p > \frac{r^a}{2}$, it follows immediately that $p$ does
not divide $a$ and if $p | r^a -1$, then $p = r^a - 1$ and $r = 2$ whence $p$ is a Mersenne prime, contradicting the
hypothesis about $p$. Therefore $p = r$ and $a = 1$. Then $|H| = p(p-1)$ and $H$ is $p$-closed by Sylow's
theorem.\qed

\begin{pro} Let $N = A_n$ and $n \ge 5$. If $N \le G = AB \le \mbox{\upshape Aut}(N)$ with
subgroups $A$ and $B$ satisfying condition (2) of the Main Theorem, $AN = BN = G$, then $A = G$ or $B = G$.
\end{pro}

{\bf Proof.} We consider first the case $n = 10$. Suppose 7 divides $|A|$. By \cite{Atl}, a 5-Sylow and a 7-Sylow subgoup of $N = A_{10}$
generate $N$. Hence if $A$ contains a 5-Sylow subgroup of $N$, then $N \le A$ whence $A = G$. So we may suppose that
$5 | |B|$. Again by \cite{Atl}, 5 and 7 are independent, contradicting the hypothesis of the proposition.

Now let $n \ne 10$. By Lemma \ref{prime}, there exists a non-Mersenne prime $p$ with $\frac{n}{2} < p \le n$. We may
assume that $p$ divides $|A|$. Take $P$ a subgroup of order $p$ in $A$. If also $p$ divides $|B|$, we
may assume by Lemma \ref{con} that $P \le B$ by replacing $B$ by a suitable conjugate if necessary.\\
Let $b \in B$ be a $q$-element for a prime $q \ne p$. By hypothesis, $\langle P, b \rangle$ is soluble.
Therefore Lemma \ref{pclo} implies that $\langle P,b \rangle$ is $p$-closed. It follows that $P$ is normalized
by all $p'$-elements of $B$ and hence by $O^{p}(B)$. If $p$ does not divide $|B|$, $B = O^p(B)$; otherwise
$B = O^p(B)P$. Therefore $B$ normalizes $P$. Hence $P \le A^b$ for all $b \in B$. It
follows now from $G = AB$ that $P \le \bigcap\limits_{b \in B}A^b = \bigcap\limits_{g \in G}A^g$. But then Core$_G(A) \ne 1$
whence $N \le A$, $A = G$.\qed

\subsection{Classical groups of Lie type}\label{class}

According to  our general strategies we proceed here as follows:

\begin{itemize}
\item[-] As mentioned in (S1), in the treatment of classical almost simple groups of Lie type, apart from two cases (see Section~2.2.2, case m); Section~2.2.3, case g)) we will always determine two primes,
one dividing $|A \cap N|$ and the other dividing $|B \cap N|$ that are independent with respect to the simple group $N$
(cf. Theorem \ref{fact}).

\item[-] Also as said in  (S4), for the infinite families of factorizations presented in \cite{LPS}, one of the primes is usually a primitive
prime divisor $r$ of $p^k - 1$ for a suitable $k$, depending on the parameters of the simple groups of Lie type of characteristic $p$. In this case
the independence of $r$ and the other specified prime, say $s$, can in general be proved by using results of \cite{ACK} where
necessary conditions for primes dividing the order of a soluble subgroup containing an element of order $r$ are given;
we will state them explicitely at the appropriate places but omit the easy calculations needed to show that $s$ is not
among the possible primes.

\item[-] There are also exceptional factorizations of classical groups with certain small parameters $n,q$.
As mentioned before in (S2), in these cases the independence of the specified primes can be deduced from the subgroup structure given in
\cite{BHR} or \cite{Atl}.\smallskip

\item[-] We also note that many of the independency results can alternatively be inferred from the proof of Theorem 2 in \cite{AI} or
from \cite{I}.
\end{itemize}

In the following we use the notation of \cite{LPS} for the groups of Lie type, in particular for the orthogonal groups.

\subsubsection{Linear groups}\label{lin}

$N = L_n(q)$, $n \ge 2$, $q = p^e$, $p$ prime.\\

We will proceed by considering the possible maximal factorizations according to \cite{LPS}, Tables 1 and 3. Using the notation of \cite{LPS},
in general $X_{\widetilde{A}} = \widetilde{A} \cap N$ and likewise for $\widetilde{B}$ as can be seen from chapters 2 and 3 of
\cite{LPS}; in the exceptional cases, we give information about $|(\widetilde{A} \cap N):X_{\widetilde{A}}|$ or
$|(\widetilde{B} \cap N):X_{\widetilde{B}}|$, respectively. Then:

\begin{itemize}\item[-] In all cases we determine independent primes, one dividing $|A \cap N|$, the other one dividing $|B \cap N|$.

\item[-] The independency of the given primes can be proved as explained at the beginning of Section \ref{class}.

\item[-] In addition, in case of $N = L_2(q)$ the independence of the given
primes is always a consequence of Dickson's list of subgroups of $L_2(q)$ as presented
for instance in Theorem II.8.27 of \cite{Hu}.
\end{itemize}

In order to deal with infinite families of factorizations we define primes $r, s, t$, depending on $(n,p,e)$, in the following way:

\begin{primes}\label{rstlg} Primes $r,\ s,\ t$ for linear groups: \begin{itemize}\item  Let $(n,p,e) \ne (2, \mbox{Mersenne prime}, 1), (2,2,3), (3,2,2), (6,2,1)$.

By Lemma \ref{Zs} there exists a primitive
prime divisor $r$ of $p^{en}-1$.

\item Let $(n,p,e) \ne (2,p,1), (2, \mbox{Mersenne prime}, 2), (2,2,6),  (4,2,2), (7,2,1),$ $(3, \mbox{Mersenne prime}, 1), (3,2,3)$.

    By Lemma \ref{Zs} there exists a primitive
prime divisor $s$ of $p^{e(n-1)}-1$.\smallskip

\item Let $(n,p,e) \ne (4, \mbox{Mersenne prime}, 1), (4,2,3), (5,2,2), (8,2,1)$, $n \ge 4$.

By Lemma \ref{Zs} there exists a primitive
prime divisor $t$ of $p^{e(n-2)}-1$.
\end{itemize}
\end{primes}

It follows from Lemma \ref{Zs} that under the given conditions $r, s, t$ do not divide $|\mbox{Out}(N)|$.\medskip

For the mentioned use of \cite{ACK}, the result that
is needed is the following:

\begin{lem}\textup{(\cite[Lemma~2.5]{ACK})}\label{rsilg} With the previous notation, let the prime $r$ be as above. If $H$ is a maximal soluble subgroup of $GL_n(q)$, $n \ge 2$, whose order is divisible by $r$, then one of the
following holds:\\
(1) $\pi(H) = \pi(n) \cup \pi(q^n -1)$;\\
(2) $\pi(H) \subseteq \pi(q-1)\cup \pi(l) \cup \{2,r\}$, where $n = 2^l$, $r = n+1$ and $q = p$ is a prime.\end{lem}

This result is sufficient for the subsequent treatment, but for some cases Lemma~2.6 of \cite{ACK} can be used alternatively.\medskip

We consider now the possible maximal factorizations as mentioned above.

\begin{itemize}
\item[{\bf a)}] $\widetilde{A} \cap N = \hat{}GL_a(q^b).b$, $ab = n$, $b$
prime, $\widetilde{B} \cap N = P_1$ or $P_{n-1}$, $(n,q) \ne (4,2)$: \smallskip \\
Here it turns out that there are always independent primes (with respect to $N$) dividing $|A \cap N|$ and $|B \cap N|$,
respectively. In Table~\ref{tlg1} we list the corresponding primes (for $|A \cap N|$ in the third column and
for $|B \cap N|$ in the fourth). Here $r$, $s$, $t$ are the primes in Notation~\ref{rstlg}, unless otherwise specified. There are eight subcases depending on the parameters $(n,p,e)$ (second column).

\begin{table}[H]
\centering
\[ \begin{array}{|l|l|c|c|} \hline
\mbox{{\bf a1)}} & n \ge 3, (n,p,e) \ne (3,\mbox{Mersenne prime},1), & r & s \\
& (3,2,3), (7,2,1),\, e > 1 \mbox{ or } s\ne n& & \\
& (r \mbox{ as above or } r = 7 \mbox{ if } (n,p,e) = (3,2,2) & &\\
& \mbox{ or }(6,2,1), & &\\
& \, \, s \mbox{ as above or } s = 7 \mbox{ if } (n,p,e) = (4,2,2)) & &\\ \hline
\mbox{{\bf a2)}} & n \ge 4,  e = 1, n = s & r & t \\
& (s \mbox{ as above or }s = 7 \mbox{ if } (n,p,e) = (7,2,1))& & \\ \hline
\mbox{{\bf a3)}} & (n,p,e) = (3,p,1), p \ne 3 & r & p \\ \hline
\mbox{{\bf a4)}} & (n,p,e) = (3,3,1) & 13 & 2 \\ \hline
\mbox{{\bf a5)}} & (n,p,e) = (3,2,3) & 73 & 7 \\ \hline
\mbox{{\bf a6)}} & (n,p,e) = (2,p,e), e \ge 2, & r & s\\
&  (n,p,e) \ne (2, \mbox{Mersenne prime}, 2)& (3 \mbox{ if } & (7 \mbox{ if }  \\
& & p = 2 \mbox{ and } e = 3) & p = 2 \mbox{ and } e = 6) \\ \hline
\mbox{{\bf a7)}} & (n,p,e) = (2, \mbox{Mersenne prime}, 2), p > 3 & r & p \\ \hline
\mbox{{\bf a8)}} & (n,p,e) = (2,p,1), p > 3 & \mbox{ largest prime} & p \\
& & \mbox{ divisor of } p + 1 &  \\ \hline
\end{array}
\]\vspace{-.4cm}
\caption{Table 1}\label{tlg1}
\end{table}

We show details for the proof of the case {\bf a1)}, as model of the arguments and strategies used along the present proof of Theorem~\ref{as}.

{\bf Case a1)} We aim to prove under the conditions of {\bf a1)} that  $r$ divides $|A\cap N|$, $s$ divides $|B\cap N|$, and $r$ and $s$ are independent primes with respect to $N$.

We first notice that, attending Notation~\ref{rstlg},  in the present case there exists a primitive prime divisor $r$ of $p^{en}-1$, except when $(n,p,e) = (3,2,2)\ \mbox{ or } (6,2,1)$, and there exists a primitive prime divisor $s$ of $p^{e(n-1)}-1$, except when $(n,p,e) = (4,2,2)$, in which case we take $s=7$.

Assume first that $(n,p,e)\neq (3,2,2),\ (6,2,1)$, and take primes $r,\ s$ as mentioned. Then we can use Lemma~\ref{rsilg} to prove the independency of $r$ and $s$.

We claim that $s\notin \pi(n)\cup \pi(q^n-1)$, and $s\notin \pi(q-1)\cup \pi(l)\cup \{2,r\}$ when $n=2^l$, $r=n+1$, $q=p$. This will imply that $r$ and $s$ are independent by Lemma~\ref{rsilg}.

This is easily checked if $(n,p,e) = (4,2,2)$ and $s=7$. Otherwise, $s$ is a primitive prime divisor of $p^{e(n-1)}-1$. By Lemma~\ref{Zs} we have that $s\ge e(n-1)+1$. If $s\in \pi(n)$, then $e=1$ and $s=n$, which is not the case. Then $s\notin \pi(n)$. On the other hand, if $s\in \pi(q^n-1)$, then $s$ divides $(q^{n-1}-1,\, q^n-1)=q-1$, which is not the case by the definition of $s$ and the fact that $n-1\ge 2$. Thus $s\notin \pi(q^n-1)\supseteq\pi(q-1)$. It is also clear that $s\neq r,\ 2$. Finally, if $n=2^l$, since $s\ge n> l$, we have that $s\notin \pi(l)$, and the claim is proven.
\smallskip

We use Lemma~\ref{l1} to prove that $r$ divides $|A\cap N|$ and $s$ divides $|B\cap N|$.

In the present case, one checks that $r$ does not divide $|\widetilde{B} \cap N||\mbox{\upshape Out}(N)|$ but divides $|N|$, which implies that $r$ divides $|A\cap N|$ by Lemma~\ref{l1}. Analogously, $s$ divides $|B\cap N|$ as it does not divide $|\widetilde{A} \cap N||\mbox{\upshape Out}(N)|$ and divides $|N|$.

Finally we consider the cases $(n,p,e) = (3,2,2)\ \mbox{ or } (6,2,1)$.

Here we take $r=7$ and notice that for the case $(n,p,e) = (3,2,2)$, $s=5$ is  primitive prime divisor of $p^{e(n-1)}-1=2^4-1$, whereas for the case $(n,p,e) = (6,2,1)$, $s=31$ is primitive prime divisor of $p^{e(n-1)}-1=2^5-1$.

In the case $(n,p,e)=(6,2,1)$, $7^2$ does not divide $|\widetilde{B} \cap N|$ but divides $|N|$, and $7$ does not divide $|\mbox{\upshape Out}(N)|$. Then Lemma~\ref{l1} implies that $7$ divides $|A\cap N|$.

In the case $(n,p,e) = (3,2,2)$, $7$ does not divide $|\widetilde{B} \cap N||\mbox{\upshape Out}(N)|$ but divides $|N|$, which implies again that $7$ divides $|A\cap N|$.

Regarding $s$, in the both cases $(n,p,e) = (3,2,2)\ \mbox{ or } (6,2,1)$, $s$ does not divide $|\widetilde{A} \cap N||\mbox{\upshape Out}(N)|$ but divides $|N|$, and Lemma~\ref{l1} implies again that $s$ divides $|B\cap N|$.

In order to prove the independency of $r$ and $s$ in the cases $(n,p,e) = (3,2,2)\ \mbox{ or } (6,2,1)$, we notice that a soluble subgroup of $N$ whose order is divisible by $rs$, would contain an $\{r,s\}$-subgroup with order also divisible by $rs$. It is easily checked, for instance using \cite{Atl}, that $N$ contains no  such $\{r,s\}$-subgroup.

\item[{\bf b)}] $X_{\widetilde{A}} = PSp_n(q)$, $\widetilde{B} \cap N = P_1$ or $P_{n-1}$, $n$ even, $n \ge 4$:\\
(Note that $|(\widetilde{A} \cap N):X_{\widetilde{A}}| = 1$ or $2$ since $PSp_n(q)$ or $PSp_n(q).2$ is a maximal
subgroup of $N$, depending on $n$ and $q$; see \cite{BHR} and \cite{KL}.) \smallskip \\
This case is handled exactly as {\bf a1)}. All other cases in {\bf a)} do not occur here since
$n \ge 4$ is even.

\item[{\bf c)}] $X_{\widetilde{A}} = PSp_n(q)$, $\widetilde{B} \cap N = Stab(V_1 \oplus V_{n-1})$, $n$ even, $n \ge 4$: \smallskip \\
This case is handled exactly as {\bf a1)}.

\item[{\bf d)}] $\widetilde{A} \cap N = \hat{}GL_{n/2}(q^2).2$, $\widetilde{B} \cap N = Stab(V_1 \oplus V_{n-1})$,
$q = 2,4$, $n$ even, $n \ge 4$: \smallskip \\
This case is handled exactly as {\bf a1)}. \medskip

We treat now all exceptional factorizations of $N$ (Table 3 of \cite{LPS}) in Table~\ref{tlg2} where the entries have the same meaning as in Table~\ref{tlg1} of {\bf a)}. In the second column we also include
$\widetilde{A} \cap N$ and $\widetilde{B} \cap N$.
\begin{table}[H]
\centering
\[ \begin{array}{|l|l|c|c|} \hline
{\rule[0mm]{0mm}{4.5mm}\mbox{{\bf e)}}} & N = L_2(11), \widetilde{A} \cap N = P_1, \widetilde{B} \cap N = A_5 & 11 & 3 \\ \hline
{\rule[0mm]{0mm}{4.5mm}\mbox{{\bf f)}}} & N = L_2(19), \widetilde{A} \cap N = P_1, \widetilde{B} \cap N = A_5 & 11 & 5 \\ \hline
{\rule[0mm]{0mm}{4.5mm}\mbox{{\bf g)}}} & N = L_2(29), \widetilde{A} \cap N = P_1, \widetilde{B} \cap N = A_5 & 29 & 5 \\ \hline
{\rule[0mm]{0mm}{4.5mm}\mbox{{\bf h)}}} & N = L_2(59), \widetilde{A} \cap N = P_1, \widetilde{B} \cap N = A_5 & 59 & 5 \\ \hline
{\rule[0mm]{0mm}{4.5mm}\mbox{{\bf i)}}} & N = L_2(7), \widetilde{A} \cap N = P_1, \widetilde{B} \cap N = S_4 & 7 & 2 \\ \hline
{\rule[0mm]{0mm}{4.5mm}\mbox{{\bf j)}}} & N = L_2(23), \widetilde{A} \cap N = P_1, \widetilde{B} \cap N = S_4 & 23 & 3 \\ \hline
{\rule[0mm]{0mm}{4.5mm}\mbox{{\bf k)}}} & N = L_2(11), \widetilde{A} \cap N = P_1, \widetilde{B} \cap N = A_4 & 11 & 3 \\ \hline
{\rule[0mm]{0mm}{4.5mm}\mbox{{\bf l)}}} & N = L_2(16), \widetilde{A} \cap N = D_{34}, \widetilde{B} \cap N = L_2(4) & 17 & 5 \\ \hline
{\rule[0mm]{0mm}{4.5mm}\mbox{{\bf m)}}} & N = L_3(4), \widetilde{A} \cap N = L_2(7), \widetilde{B} \cap N = A_6 & 7 & 5 \\ \hline
{\rule[0mm]{0mm}{4.5mm}\mbox{{\bf n)}}} & N = L_5(2), \widetilde{A} \cap N = 31.5, \widetilde{B} \cap N = P_2 \mbox{ or } P_3 & 31 & 7 \\ \hline
\end{array}
\]\vspace{-.4cm}
\caption{Table 2}\label{tlg2}
\end{table}

{\bf Notation for tables.}  In all tables the meaning of the entries is the same as in Table~\ref{tlg1} of {\bf a)}, i.e. the entries in the third and fourth columns are independent primes (with respect to $N$)  dividing $|A \cap N|$ (third column) and $|B \cap N|$ (fourth column), respectively.
\medskip

\subsubsection{Unitary groups}\label{un}

$N = U_n(q)$, $n \ge 3$, $q = p^e$, $p$ prime.\\

Infinite families of factorizations of $N$ described in Table 1 of \cite{LPS} exist only for even $n$. So we assume first that
$n = 2m \ge 4$. Then:

\begin{itemize}\item[-] In all cases we determine independent primes, one dividing $|A \cap N|$,
the other one dividing $|B \cap N|$.
\item[-] The divisibility properties and the independency are proved as described at the beginning of Section \ref{class}.
\end{itemize}

We define corresponding convenient primes $r, s$, depending on $(n,p,e)$, as follows:

\begin{primes}\label{rsug} Primes $r,\ s$ for unitary groups: \begin{itemize}
\item[$\bullet$] Let $(n,p,e) \ne (4,2,1)$.
By Lemma \ref{Zs} there exists a primitive
prime divisor $r$ of $p^{2e(n-1)}-1$.

\item[$\bullet$] Let $(n,p,e) \ne (6,2,1)$.
By Lemma \ref{Zs} there exists a primitive
prime divisor $s$ of $p^{en}-1$. \end{itemize}\end{primes}

It follows from Lemma \ref{Zs} that under the given conditions $r, s$ do not divide $|\mbox{Out}(N)|$.

Here we make use of the following result in \cite{ACK}:

\begin{lem}\textup{(\cite[Lemma~2.8(1)]{ACK})}\label{ACK2.8} With the previous notation, let the prime $r$ be as above. If $H$ is a maximal soluble subgroup of $U_n(q)$, $n = 2m \ge 4$, with $r \, | \, |H|$, then
$\pi(H) \subseteq \pi(n-1) \cup \pi(q^{n-1} +1)$.\end{lem}

We consider now the possible maximal factorizations as mentioned before.

\begin{itemize}
\item[{\bf a)}] $\widetilde{A} \cap N = N_1$, $\widetilde{B} \cap N = P_m$: \smallskip\\
We have to consider two subcases, presented in Table~\ref{tug1}.
\begin{table}[H]
\centering
\[ \begin{array}{|l|l|c|c|} \hline
\mbox{{\bf a1)}} & (n,p,e) \ne (4,2,1), & r & s \, (s \mbox{ as above or } s = 7 \mbox{ if } (n,p,e) = (6,2,1))\\ \hline
\mbox{{\bf a2)}} & (n,p,e) = (4,2,1) & 3 & 5 \\ \hline
\end{array}
\]\vspace{-.4cm}
\caption{Table 3}\label{tug1}
\end{table}

This case is proven by using analogous arguments to those for linear groups, in Section~2.2.1, case {\bf a1)}, but with corresponding primes $r, \ s$ as in Notation~\ref{rsug} or as specified, and Lemma~\ref{ACK2.8} for the proof of the independence.

\item[{\bf b)}] $\widetilde{A} \cap N = N_1$, $X_{\widetilde{B}} = PSp_n(q)$:\\
Note that $|(\widetilde{B} \cap N): X_{\widetilde{B}}| = 1$ or $2$ since $PSp_n(q)$ or $PSp_n(q).2$ is a maximal
subgroup of $N$ (depending on $q$ and $n$; cf. \cite{KL} and \cite{BHR}). \smallskip \\
This case is handled exactly as {\bf a)}.

\item[{\bf c)}] $N = U_n(2)$, $m \ge 3$, $\widetilde{A} \cap N = N_1$, $X_{\widetilde{B}} = \hat{}SL_m(4).2$: \smallskip \\
This case is handled exactly as {\bf a1)}.

\item[{\bf d)}] $N = U_n(4)$, $\widetilde{A} \cap N = N_1$, $X_{\widetilde{B}} = \hat{}SL_m(16).3.2$: \smallskip \\
This case is handled exactly as {\bf a1)}.

\end{itemize}

We treat now the exceptional factorizations as given in Table 3 of \cite{LPS}.\smallskip

 Here we have:
\begin{itemize}\item[-] Apart from one case ($N=U_4(3)$) there are independent
primes dividing $|A \cap N|$ and $|B \cap N|$, respectively, summarized in Table~\ref{tug2}.
\item[-] The divisibility of $|A \cap N|$ and $|B \cap N|$
by the specified primes follows from Lemma \ref{l1}.\end{itemize}
\begin{table}[H]
\centering
\[ \begin{array}{|l|l|c|c|} \hline
{\rule[0mm]{0mm}{4.5mm}\mbox{{\bf e)}}} & N = U_3(3), \widetilde{A} \cap N = L_2(7), \widetilde{B} \cap N = P_1 & 7 & 2 \\ \hline
{\rule[0mm]{0mm}{4.5mm}\mbox{{\bf f)}}} & N = U_3(5), \widetilde{A} \cap N = A_7, \widetilde{B} \cap N = P_1 & 7 & 5 \\ \hline
{\rule[0mm]{0mm}{4.5mm}\mbox{{\bf g)}}} & N = U_3(8), \widetilde{A} \cap N = 19.3, \widetilde{B} \cap N = P_1 & 19 & 7 \\ \hline
{\rule[0mm]{0mm}{4.5mm}\mbox{{\bf h)}}} & N = U_4(2), \widetilde{A} \cap N = 3^3.S_4, \widetilde{B} \cap N = P_2 & 3 & 5 \\ \hline
{\rule[0mm]{0mm}{4.5mm}\mbox{{\bf i)}}} & N = U_6(2), \widetilde{A} \cap N = N_1, \widetilde{B} \cap N = U_4(3).2 & 11 & 7 \\ \hline
{\rule[0mm]{0mm}{4.5mm}\mbox{{\bf j)}}} & N = U_6(2), \widetilde{A} \cap N = N_1, \widetilde{B} \cap N = M_{22} & 11 & 7 \\ \hline
{\rule[0mm]{0mm}{4.5mm}\mbox{{\bf k)}}} & N = U_9(2), \widetilde{A} \cap N = J_3, \widetilde{B} \cap N = P_1 & 19 & 43 \\ \hline
{\rule[0mm]{0mm}{4.5mm}\mbox{{\bf l)}}} & N = U_{12}(2), \widetilde{A} \cap N = Suz, \widetilde{B} \cap N = N_1 & 13 & 683 \\ \hline
\end{array}
\]\vspace{-.4cm}
\caption{Table 4}\label{tug2}
\end{table}
For case j) we have to explain why 11 divides $|A \cap N|$. If not, this would lead to a non-trivial factorization
$A(\widetilde{A} \cap \widetilde{B})$ of $\widetilde{A}$, $\widetilde{A}$ an almost simple group with socle
$N_1 \cong U_5(2)$. This is impossible by \cite{LPS}.\medskip

It remains to consider

\begin{itemize}
\item[{\bf m)}] $N = U_4(3)$, $\widetilde{A} \cap N = L_3(4)$, $\widetilde{B} \cap N = P_1$ or $PSp_4(3)$ or $P_2$: \smallskip \\
It follows from Lemma \ref{l1} that in all three cases 7 divides $|A \cap N|$
and $3^4$ divides $|B \cap N|$. \\
Suppose that $A$ and $B$ satisfy condition (2) of the Main Theorem. Note that 7 is a primitive prime divisor of $3^6 - 1$. Let $y \in A \cap N$
of order 7. By \cite{Atl}, $|N_N(\langle y \rangle)| = 7 \cdot 3$. This and Lemma~\ref{ACK2.8} imply that
$N_N(\langle y \rangle)$ is the only maximal soluble subgroup of $N$ containing $y$. Therefore $\langle x,y \rangle =
N_N(\langle y \rangle)$ for all non-trivial 3-elements $x$ of $B \cap N$, a contradiction.\\
(We note that in case $G = N$ the primes 2 and 7 are independent and 2 divides $|B\cap N|$. But this need not be the
case if $G \ge N.2$, see Theorem \ref{fact}.)
\end{itemize}

\subsubsection{Symplectic groups}\label{symp}

$N = PSp_{2m}(q)$, $m \ge 2$, $q = p^e$, $p$ prime.\\

For arguments like in the previous cases, we define primes $r, s, t$, depending on $(m,p,e)$, in the following way:

\begin{primes}\label{rstsg} Primes $r,\ s,\ t$ for symplectic groups: \begin{itemize}
\item[$\bullet$]  Let $(m,p,e) \ne (3,2,1)$. By Lemma \ref{Zs} there exists a primitive
prime divisor $r$ of $p^{2em}-1$.
\item[$\bullet$] Let $(m,p,e) \ne (2, \mbox{Mersenne prime}, 1), (2,2,3), (3,2,2), (6,2,1)$. By Lemma \ref{Zs} there exists a primitive
prime divisor $s$ of $p^{em}-1$.

\item[$\bullet$] Let $(m,p,e) \ne (2, \mbox{Mersenne prime}, 1), (2,2,3), (4,2,1)$. By Lemma \ref{Zs} there exists a primitive
prime divisor $t$ of $p^{2e(m-1)}-1$.
\end{itemize}\end{primes}

It follows from Lemma \ref{Zs} that under the given conditions $r, s, t$ do not divide $|\mbox{Out}(N)|$.\medskip

We first consider the possible maximal factorizations according to \cite{LPS}, Table 1. Then:

\begin{itemize}\item[-] In all but one case ($N=PSp_6(2)$) there exist independent primes with respect to $N$, one dividing $|A \cap N|$,
the other one dividing $|B \cap N|$.
\item[-]  The
independency of the given primes can be proved as mentioned at the beginning of Section \ref{class}.
\end{itemize}

The result that is needed from \cite{ACK} is the following. For later use in Sections~\ref{sso}~and~\ref{sso-} we formulate it already in a way that also includes orthogonal groups in odd dimension and of -- type in even dimension.

\begin{lem}\textup{(\cite[Lemma~2.8(2)]{ACK})}\label{lem8} With the previous notation, let the prime $r$ be as above. If $H$ is a maximal soluble subgroup of $PSp_{2m}(q)$, $m \ge 2$, $P\Omega_{2m+1}(q)$, $m\ge 3$ and $q$ odd, or $P\Omega_{2m}^{-}(q)$, $m\ge 4$, whose order is divisible by $r$, then one of the
following holds:\\
(1) $\pi(H) \subseteq \pi(m) \cup \pi(q^m +1) \cup \{2\}$;\\
(2) $\pi(H) \subseteq \pi(q-1) \cup \pi(l+1)\cup \{2,r\}$, where $m = 2^l$, $r = 2m+1$ and $q = p$ is a prime.\end{lem}

The possible maximal factorizations are as follows.

\begin{itemize}
\item[{\bf a)}] $\widetilde{A} \cap N = PSp_{2a}(q^b).b$, $ab = m$, $b$ prime, $\widetilde{B} \cap N = P_1$: \smallskip \\
In this case there are always independent primes (with respect to $N$) dividing $|A \cap N|$ and $|B \cap N|$,
respectively. They are specified in Table~\ref{tsg1}. The primes $r$ and
$t$ are as in Notation~\ref{rstsg}.
\begin{table}[H]
\centering
\[ \begin{array}{|l|l|c|c|} \hline
\mbox{{\bf a1)}} & (m,p,e) \ne (2,\mbox{Mersenne prime},1) & r & t \\
& (r = 7 \mbox{ if } (m,p,e) = (3,2,1),& & \\
& \, t = 7 \mbox{ if } (m,p,e) = (2,2,3) \mbox{ or } (4,2,1))& & \\ \hline
\mbox{{\bf a2)}} & (m,p,e) = (2,\mbox{Mersenne prime},1) & r & p \\ \hline
\end{array}
\]\vspace{-.4cm}
\caption{Table 5}\label{tsg1}
\end{table}

\item[{\bf b)}] $q = 2^e$, $\widetilde{A} \cap N = Sp_{2a}(q^b).b$, $ab = m$, $b$ prime, $\widetilde{B} \cap N = O_{2m}^+(q)$:\\
(Note that $X_{\widetilde{B}} = \widetilde{B} \cap N$ since $O_{2m}^+(q)$ is a maximal subgroup of $N = Sp_{2m}(q)$, see
\cite{KL} and \cite{BHR}.)\\

There are always independent primes (with respect to $N$) dividing $|A \cap N|$ and $|B \cap N|$,
respectively. They are specified in Table~\ref{tsg2} for all subcases to be considered. The meaning of $r$ and
$t$ is as above in Notation~\ref{rstsg}.
\begin{table}[H]
\centering
\[ \begin{array}{|l|l|c|c|} \hline
\mbox{{\bf b1)}} & (m,e) \ne (2,3), (3,1), (4,1), & r & t \\ \hline
\mbox{{\bf b2)}} & (m,e) = (2,3) & 13 & 7 \\ \hline
\mbox{{\bf b3)}} & (m,e) = (3,1) & 7 & 5 \\ \hline
\mbox{{\bf b4)}} & (m,e) = (4,1) & 17 & 7 \\ \hline
\end{array}
\]\vspace{-.4cm}
\caption{Table 6}\label{tsg2}
\end{table}

In {\bf b3)} we have to justify why 7 divides $|A|$ (note that $|\mbox{Out}(N)| = 1$). Since $2^2$ and $3^2$ divide
$|A|$, this follows from the fact that there is no subgroup of $Sp_2(8).3$ whose order is divisible by 36, but not by 7.\\

\item[{\bf c)}] $q = 2^e$, $\widetilde{A} \cap N = Sp_{2a}(q^b).b$, $ab = m$, $b$ prime, $\widetilde{B} \cap N = O_{2m}^-(q)$:\\
(Note that $X_{\widetilde{B}} = \widetilde{B} \cap N$ since $O_{2m}^-(q)$ is a maximal subgroup of $N = Sp_{2m}(q)$, see
\cite{KL} and \cite{BHR}.)\\

Again there are always independent primes dividing $|A \cap N|$ and $|B \cap N|$, respectively. They are presented in Table~\ref{tsg3} where the meaning of $r$, $s$ and $t$ is as above in Notation~\ref{rstsg}.
\begin{table}[H]
\centering
\[ \begin{array}{|l|l|c|c|} \hline
\mbox{{\bf c1)}} & (m,e) \ne (2,3), (3,1), (4,1), r \mid |A \cap N| & r & t \\ \hline
\mbox{{\bf c2)}} & (m,e) \ne (2,3), (3,1), (3,2), (6,1), r \not\:\mid |A \cap N| & s & r \\ \hline
\mbox{{\bf c3)}} & (m,e) = (2,3), r = 13 \mid |A \cap N|  & 13 & 7 \\ \hline
\mbox{{\bf c4)}} & (m,e) = (2,3), r = 13 \not\:\mid |A \cap N|  & 7 & 13 \\ \hline
\mbox{{\bf c5)}} & (m,e) = (3,1) & 7 & 5 \\ \hline
\mbox{{\bf c6)}} & (m,e) = (4,1) & 17 & 7 \\ \hline
\mbox{{\bf c7)}} & (m,e) = (3,2) \mbox{ or } (6,1), r = 13 \not\:\mid |A \cap N|  & 7 & 13 \\ \hline
\end{array}
\]\vspace{-.4cm}
\caption{Table 7}\label{tsg3}
\end{table}

\item[{\bf d)}] $q = 2^e$, $\widetilde{A} \cap N = O_{2m}^-(q)$, $\widetilde{B} \cap N = P_m$: \smallskip \\
Here there are always independent primes dividing $|A \cap N|$ and $|B \cap N|$, respectively. They are presented in Table~\ref{tsg4} where the meaning of $r$ and $s$ is as above in Notation~\ref{rstsg}
\begin{table}[H]
\centering
\[ \begin{array}{|l|l|c|c|} \hline
\mbox{{\bf d1)}} & (m,e) \ne (2,3), (3,1) & r & s \\
& (s = 7 \mbox{ for } (m,e) = (3,2) \mbox{ or } (6,1)) & & \\ \hline
\mbox{{\bf d2)}} & (m,e) = (2,3) & 13 & 7 \\ \hline
\mbox{{\bf d3)}} & (m,e) = (3,1) & 5 & 7 \\ \hline
\end{array}
\]\vspace{-.4cm}
\caption{Table 8}\label{tsg4}
\end{table}
\item[{\bf e)}] $q = 2^e$, $m$ even, $\widetilde{A} \cap N = O_{2m}^-(q)$, $\widetilde{B} \cap N = Sp_m(q)\, wr \, S_2$: \smallskip \\
This case case is handled exactly like {\bf d)}.
\item[{\bf f)}] $q = 2$, $m$ even, $\widetilde{A} = Sp_m(2).2$, $\widetilde{B} = N_2$ (note that $|\mbox{Out}(N)| = 1$): \smallskip\\
One can assume that $m \ge 4$ since $Sp_4(2)' = L_2(9)$. Then $r$ divides $|A \cap N|$ and $t$ divides $|B \cap N|$
where $r$ and $t$ are as in Notation~\ref{rstsg}, or $t = 7$ if $m = 4$. $r$ and $t$ are independent.
\item[{\bf g)}] $q = 2$, $\widetilde{A} = O_{2m}^-(2)$, $\widetilde{B} = O_{2m}^+(2)$ (note that $|\mbox{Out}(N)| = 1$): \smallskip \\
If $m \ne 3$, let $r$ and $s$ as above in Notation~\ref{rstsg}, or $s = 7$ for $m =6$. Then $r$ divides $|A|$ and $s$ divides $|B|$.
$r$ and $s$ are independent.

Let $m = 3$. By Lemma \ref{l1}, $3^2$ divides $|A|$ and 7 divides $|B|$. If 5 divides $|A|$,
we are done since 5 and 7 are independent. Otherwise 5 divides $|B|$. Using \cite{BHR} or \cite{Atl},
the only maximal soluble subgroups of $N$ whose order is divisible by 15 are the normalizers of
elements of order 5, isomorphic to $(5 \times S_3):4$. Suppose $A$ and $B$ satisfy condition (2) of the Main Theorem.
It follows that if $y \in B$ is of order 5, then $\langle x_1,y \rangle = \langle x_2,y \rangle$ of order 15 for all non-trivial 3-elements
$x_1,x_2 \in A$, a contradiction.\\
(We note that in case $m = 3$ there need not be independent primes dividing $|A|$ and $|B|$, respectively. See
Theorem \ref{fact}.)
\end{itemize}
It remains to consider the three cases {\bf h)}, {\bf i)}, {\bf j)} summarized in Table~\ref{tsg5}.

Note in case {\bf j)} that $X_{\widetilde{B}} = \widetilde{B} \cap N$ by \cite{LPS}, 3.2.4(d).

\begin{table}[H]
\centering
\[ \begin{array}{|l|l|c|c|} \hline
{\rule[0mm]{0mm}{4.5mm}\mbox{{\bf h)}}} & q = 4, \widetilde{A} \cap N = O_{2m}^-(4), \widetilde{B} \cap N = O_{2m}^+(4) & r & s \,(s = 7 \mbox{ if } m = 3) \\ \hline
{\rule[0mm]{0mm}{4.5mm}\mbox{{\bf i)}}} & q = 4, m \mbox{ even, } \widetilde{A} \cap N = Sp_m(16).2, \widetilde{B} \cap N = N_2 & r & t \\ \hline
{\rule[0mm]{0mm}{4.5mm}\mbox{{\bf j)}}} & q = 4, 16,  \widetilde{A} \cap N = O_{2m}^-(q), \widetilde{B} \cap N = Sp_{2m}(q^{1/2}) & r & s \\ \hline
\end{array}
\] \vspace{-.4cm}
\caption{Table 9}\label{tsg5}
\end{table}

We now consider the possible maximal factorizations according to \cite{LPS}, Table 2. Here:
\begin{itemize}\item[-] There are always
independent primes dividing $|A \cap N|$ and $|B \cap N|$, respectively, as specified next.
\end{itemize}

\begin{itemize}
\item[{\bf k)}] $m = 2$, $q = 2^e$, $e$ odd, $e \ge 3$, $\widetilde{A} \cap N = Sz(q)$, $\widetilde{B} \cap N = O_4^+(q)$: \smallskip \\
$r$ divides $|A \cap N|$, $s$ divides $|B \cap N|$, $r$ and $s$ are independent.

\item[{\bf l)}] $m = 3$, $q = 2^e$, $\widetilde{A} \cap N = G_2(q)$, $\widetilde{B} \cap N = O_6^+(q)$: \smallskip \\
If $e > 1$, then $r$ divides $|A \cap N|$, $t$ divides $|B \cap N|$, $r$ and $t$ are independent.\\
If $e = 1$, then $G = N$. Clearly 5 divides $|B|$. Since $N = A \widetilde{B}$, it follows that
$G_2(2) \cong \widetilde{A} = A(\widetilde{A} \cap \widetilde{B})$. But according to \cite{LPS}, Table 5, $G_2(2)$ has
only trivial factorizations, whence $A = \widetilde{A}$ and 7 divides $|A|$. 7 and 5 are independent.

\item[{\bf m)}] $m = 3$, $q = 2^e$, $\widetilde{A} \cap N = G_2(q)$, $\widetilde{B} \cap N = O_6^-(q)$: \smallskip \\
Assume first the $e > 1$. $t$ divides $|B \cap N|$. If $r$ divides $|A \cap N|$, we choose the independent primes $r$
and $t$. If $r$ does not divide $|A \cap N|$, then $r$ divides $|B \cap N|$. As $s$ divides $|A \cap N|$ ($s = 7$
if $e = 2$), we choose the independent primes $s$ and $r$.\\
If $e = 1$, 7 divides $|A \cap N|$ and 5 divides $|B \cap N|$. 7 and 5 are independent.

\item[{\bf n)}] $m = 3$, $q = 2^e$, $\widetilde{A} \cap N = G_2(q)$, $\widetilde{B} \cap N = P_1$ or $N_2$: \smallskip \\
This follows exactly like case {\bf l)}.

\end{itemize}

We now consider the exceptional maximal factorizations according to \cite{LPS}, Table 3. Here:
\begin{itemize}\item[-] There are always
independent primes dividing $|A \cap N|$ and $|B \cap N|$, respectively, specified in Table~\ref{tsg6}.\end{itemize}
Note that Out$(N) = 1$ in {\bf q)} and {\bf r)}.
\begin{table}[H]
\centering
\[ \begin{array}{|l|l|c|c|} \hline
{\rule[0mm]{0mm}{4.5mm}\mbox{{\bf o)}}} & N = PSp_4(3), \widetilde{A} \cap N = 2^4.A_5, \widetilde{B} \cap N = P_1 \mbox{ or } P_2 & 5 & 3 \\ \hline
{\rule[0mm]{0mm}{4.5mm}\mbox{{\bf p)}}} & N = PSp_6(3), \widetilde{A} \cap N = L_2(13), \widetilde{B} \cap N = P_1 & 7 & 5 \\ \hline
{\rule[0mm]{0mm}{4.5mm}\mbox{{\bf q)}}} & N = Sp_8(2), \widetilde{A} = O_8^-(2), \widetilde{B} = S_{10} & 17 & 5 \\ \hline
{\rule[0mm]{0mm}{4.5mm}\mbox{{\bf r)}}} & N = Sp_8(2), \widetilde{A} = L_2(17), \widetilde{B} = O_8^+(2) & 17 & 5 \\ \hline
\end{array}
\] \vspace{-.4cm}
\caption{Table 10}\label{tsg6}
\end{table}

\subsubsection{Orthogonal groups in odd dimension}\label{sso}

$N = P{\Omega}_{2m+1}(q) \, (= \Omega_{2m+1}(q))$, $m \ge 3$, $q = p^e$ odd, $p$ prime.\\

We consider the maximal factorizations according to \cite{LPS}, Tables 1,2,3. Then:\smallskip

\begin{itemize}
 \item[-] There are always
independent primes dividing $|A \cap N|$ and $|B \cap N|$, respectively, specified in Table~\ref{tog}.
\end{itemize}

In order to prove it, we consider primes as follows, and   Lemma~\ref{lem8}.

\begin{primes}\label{rstog} Primes $r,\ s,\ t$ for orthogonal groups in odd dimension:  \begin{itemize}\item[$\bullet$]
Let $r$, $s$ and $t$ be primes as in Notation~\ref{rstsg}.\end{itemize}\end{primes}

Note that they exist in this case without any restriction since $m \ge 3$ and $p \ne 2$.
Note that none of $r$, $s$, $t$ divides $|\mbox{Out}(N)|$.
\begin{table}[H]
\centering
\[ \begin{array}{|l|l|c|c|} \hline
{\rule[0mm]{0mm}{4.5mm}\mbox{{\bf a)}}} &  \widetilde{A} \cap N = N_1^-, \widetilde{B} \cap N = P_m & r & s \\ \hline
{\rule[0mm]{0mm}{4.5mm}\mbox{{\bf b)}}} & N = \Omega_7(q), \widetilde{A} \cap N = G_2(q), \widetilde{B} \cap N = P_1 & r & t \\ \hline
{\rule[0mm]{0mm}{4.5mm}\mbox{{\bf c)}}} & N = \Omega_7(q), \widetilde{A} \cap N = G_2(q), \widetilde{B} \cap N = N_1^+ & r & t \\ \hline
{\rule[0mm]{0mm}{4.5mm}\mbox{{\bf d)}}} & N = \Omega_7(q), \widetilde{A} \cap N = G_2(q), \widetilde{B} \cap N = N_1^-, r \,|\, |A \cap N| & r & t \\ \hline
{\rule[0mm]{0mm}{4.5mm}\mbox{{\bf e)}}} & N = \Omega_7(q), \widetilde{A} \cap N = G_2(q), \widetilde{B} \cap N = N_1^-, r \not\:\mid |A \cap N| & s & r \\ \hline
{\rule[0mm]{0mm}{4.5mm}\mbox{{\bf f)}}} & N = \Omega_7(q), \widetilde{A} \cap N = G_2(q), \widetilde{B} \cap N = N_2^{\epsilon}, \epsilon = +,- & r & t \\ \hline
{\rule[0mm]{0mm}{4.5mm}\mbox{{\bf g)}}} & N = \Omega_{13}(3^e), \widetilde{A} \cap N = PSp_6(3^e).a, a \le 2,  \widetilde{B} \cap N = N_1^- & r & s \\ \hline
{\rule[0mm]{0mm}{4.5mm}\mbox{{\bf h)}}} & N = \Omega_{25}(3^e), \widetilde{A} \cap N = F_4(3^e), \widetilde{B} \cap N = N_1^- & s & r \\ \hline
{\rule[0mm]{0mm}{4.5mm}\mbox{{\bf i)}}} & N = \Omega_7(3), \widetilde{A} \cap N = G_2(3), \widetilde{B} \cap N = Sp_6(2) \mbox{ or } S_9, 7 \, | \,|A \cap N| & 7 & 5 \\ \hline
{\rule[0mm]{0mm}{4.5mm}\mbox{{\bf j)}}} & N = \Omega_7(3), \widetilde{A} \cap N = G_2(3), \widetilde{B} \cap N = Sp_6(2) \mbox{ or } S_9, 7 \not\:\mid |A \cap N| & 13 & 7 \\ \hline
{\rule[0mm]{0mm}{4.5mm}\mbox{{\bf k)}}} & N = \Omega_7(3), \widetilde{A} \cap N = S_9, \widetilde{B} \cap N = N_1^+ \mbox{ or } P_3 & 7 & 13 \\ \hline
{\rule[0mm]{0mm}{4.5mm}\mbox{{\bf l)}}} & N = \Omega_7(3), \widetilde{A} \cap N = Sp_6(2), \widetilde{B} \cap N = N_1^+ \mbox{ or } P_3 & 7 & 13 \\ \hline
{\rule[0mm]{0mm}{4.5mm}\mbox{{\bf m)}}} & N = \Omega_7(3), \widetilde{A} \cap N = 2^6.A_7, \widetilde{B} \cap N = P_3 & 7 & 13 \\ \hline
\end{array}
\]\vspace{-.4cm}
\caption{Table 11}\label{tog}
\end{table}

\subsubsection{Orthogonal groups of -- type in even dimension}\label{sso-}

$N = P{\Omega}_{2m}^-(q)$, $m \ge 4$, $q = p^e$, $p$ prime.\\

We now consider the maximal factorizations according to \cite{LPS}, Tables 1, 3. Here:

\begin{itemize}\item[-] There are always
independent primes dividing $|A \cap N|$ and $|B \cap N|$, respectively, specified in Table~\ref{tog-}.
\end{itemize}

In order to prove it, we again consider primes as follows, and  Lemma~\ref{lem8}.

\begin{primes}\label{rtog-} Primes $r,\ t$ for orthogonal groups of -- type in even dimension:
 \begin{itemize}\item[$\bullet$]
Let $r$, $t$ be primes as in Notation~\ref{rstsg}.\end{itemize}\end{primes}

Note that they exist in this case with the only restriction $(m,p,e) \ne (4,2,1)$ to ensure the existence of $t$. Note that both $r$ and $t$ do not divide $|\mbox{Out}(N)|$.
\begin{table}[H]
\centering
\[ \begin{array}{|l|l|c|c|} \hline
{\rule[0mm]{0mm}{4.5mm}\mbox{{\bf a)}}} &  \widetilde{A} \cap N = P_1 \mbox{ or } N_1, \widetilde{B} \cap N = \hat{}GU_m(q), m \mbox{ odd } & t & r \\ \hline
{\rule[0mm]{0mm}{4.5mm}\mbox{{\bf b)}}} &  \widetilde{A} \cap N = N_1, \widetilde{B} \cap N = \Omega_m^-(q^2).2, q = 2,4, m \mbox{ even } & t & r \\
 &  &  (t = 7 \mbox{ for } (m,q) = (4,2)) & \\ \hline
{\rule[0mm]{0mm}{4.5mm}\mbox{{\bf c)}}} &  \widetilde{A} \cap N = N_2^+, \widetilde{B} \cap N = GU_m(4), q = 4, m \mbox{ odd } & t & r \\ \hline
{\rule[0mm]{0mm}{4.5mm}\mbox{{\bf d)}}} & N = P{\Omega}_{10}^-(2), \widetilde{A} \cap N = A_{12}, \widetilde{B} \cap N = P_1 & 11 & 17 \\ \hline
\end{array}
\]\vspace{-.4cm}
\caption{Table 12}\label{tog-}
\end{table}

\subsubsection{Orthogonal groups of + type in even dimension}

$N = P{\Omega}_{2m}^+(q)$, $m \ge 4$, $q = p^e$, $p$ prime.\\

We will consider below the maximal factorizations according to \cite{LPS}, Tables 1, 2, 3 (case $m\ge 5$) and Table 4 (case $m=4$).

For arguments like in previous cases, we define primes $r, s, t$, depending on $(m,p,e)$, in the following way:

\begin{primes}\label{rstog+} Primes $r,\ s,\ t$ for orthogonal groups of + type in even dimension:  \begin{itemize}\item[$\bullet$]
Let $(m,p,e) \ne (4,2,1)$. By Lemma \ref{Zs} there exists a primitive
prime divisor $r$ of $p^{2e(m-1)}-1$.

\item[$\bullet$] Let $(m,p,e) \ne (6,2,1)$. By Lemma \ref{Zs} there exists a primitive
prime divisor $s$ of $p^{em}-1$.

\item[$\bullet$] Let $(m,p,e) \ne (4,2,2), (7,2,1)$. By Lemma \ref{Zs} there exists a primitive
prime divisor $t$ of $p^{e(m-1)}-1$.\end{itemize}\end{primes}

It follows from Lemma \ref{Zs} that under the given conditions $r, s, t$ do not divide $|\mbox{Out}(N)|$. Moreover,
$r, s, t$ are pairwise distinct.

For independency proofs we will need the following result of \cite{ACK} (with a corrected misprint there):

\begin{lem}\textup{(\cite[Lemma~2.8(3)]{ACK})} With the previous notation, let the prime $r$ be as above. If $H$ is a maximal soluble subgroup of $N = P{\Omega}_{2m}^+(q)$, $m \ge 4$, whose order
is divisible by $r$, then one of the following holds:\\
(1) $\pi(H) \subseteq \pi(m-1) \cup \pi(q^{m-1} +1) \cup \pi(q+1) \cup \{2\}$;\\
(2) $\pi(H) \subseteq \pi(q^2-1)\cup \pi(l+1) \cup \{2,r\}$, where $m-1 = 2^l$, $r = 2m-1$ and $q = p$ is a prime;\\
(3) $\pi(H) \subseteq \pi(q-1)\cup \pi(l+1) \cup \{2,r\}$, where $m = 2^l$, $r = 2m-1$ and $q = p$ is a prime. \end{lem}

We assume first that $m \ge 5$ and  consider the maximal factorizations according to \cite{LPS}, Tables 1, 2, 3. Then:
\begin{itemize}\item[-] There are always
independent primes dividing $|A \cap N|$ and $|B \cap N|$, respectively, specified in Table~\ref{tog+}.\end{itemize}
We note in part {\bf d)} that $X_{\widetilde{B}} = PSp_2(q) \otimes PSp_m(q)$ has index 1 or 2 in $\widetilde{B} \cap N$; cf. \cite{KL}, Table~3.5.E and Prop.~4.4.12.
\begin{table}[H]
\centering
\[ \begin{array}{|l|l|c|c|} \hline
{\rule[0mm]{0mm}{4.5mm}\mbox{{\bf a)}}} &  \widetilde{A} \cap N = N_1, \widetilde{B} \cap N = P_m \mbox{ or } P_{m-1} & r & s \\ \hline
{\rule[0mm]{0mm}{4.5mm}\mbox{{\bf b)}}} &  \widetilde{A} \cap N = N_1, \widetilde{B} \cap N = \hat{}GU_m(q).2, & r & s \\
 & m \mbox{ even, } r \,| \,|A \cap N| & & \\ \hline
{\rule[0mm]{0mm}{4.5mm}\mbox{{\bf c)}}} &  \widetilde{A} \cap N = N_1, \widetilde{B} \cap N = \hat{}GU_m(q).2, & t & r \\
 &  m \mbox{ even, } r \not\:\mid |A \cap N| & & \\ \hline
{\rule[0mm]{0mm}{4.5mm}\mbox{{\bf d)}}} &  \widetilde{A} \cap N = N_1, X_{\widetilde{B}} = PSp_2(q) \otimes PSp_m(q), & r & s \\
 & m \mbox{ even, } q > 2 & & \\ \hline
{\rule[0mm]{0mm}{4.5mm}\mbox{{\bf e)}}} &  \widetilde{A} \cap N = N_2^-, \widetilde{B} \cap N = P_m \mbox{ or } P_{m-1} & r & s \\
 &  &  &(s = 7 \mbox{ for } (m,q) = (6,2)) \\ \hline
{\rule[0mm]{0mm}{4.5mm}\mbox{{\bf f)}}} &  \widetilde{A} \cap N = P_1, \widetilde{B} \cap N = \hat{}GU_m(q).2, m \mbox{ even } & t & r \\ \hline
{\rule[0mm]{0mm}{4.5mm}\mbox{{\bf g)}}} &  \widetilde{A} \cap N = N_1, \widetilde{B} \cap N = \hat{}GL_m(q).2 & r & s \\
 &  &  &(s = 7 \mbox{ for } (m,q) = (6,2)) \\ \hline
 \end{array}
\]
\end{table}
\begin{center}\vspace{-.5cm} Table 13 (contd.)\end{center}

\begin{table}[H]
\centering
\[ \begin{array}{|l|l|c|c|} \hline
{\rule[0mm]{0mm}{4.5mm}\mbox{{\bf h)}}} &  \widetilde{A} \cap N = N_1,  \widetilde{B} \cap N = \Omega_m^+(4).2^2, q = 2, & r & s \\
 &  m \mbox{ even } &  &(s = 7 \mbox{ for } (m,q) = (6,2)) \\ \hline
{\rule[0mm]{0mm}{4.5mm}\mbox{{\bf i)}}} &  \widetilde{A} \cap N = N_1,  \widetilde{B} \cap N = \Omega_m^+(16).2^2, q = 4, & r & s \\
 &  m \mbox{ even } & & \\ \hline
{\rule[0mm]{0mm}{4.5mm}\mbox{{\bf j)}}} &  \widetilde{A} \cap N = N_2^-, \widetilde{B} \cap N = \hat{}GL_m(2).2, q = 2 & r & s \\
 &  &  &(s = 7 \mbox{ for } (m,q) = (6,2))\\ \hline
{\rule[0mm]{0mm}{4.5mm}\mbox{{\bf k)}}} &  \widetilde{A} \cap N = N_2^-, \widetilde{B} \cap N = \hat{}GL_m(4).2, q = 4 & r & s \\ \hline
{\rule[0mm]{0mm}{4.5mm}\mbox{{\bf l)}}} &  \widetilde{A} \cap N = N_2^+, \widetilde{B} \cap N = \hat{}GU_m(4).2, q = 4, & t & r \\
 &  m \mbox{ even } & & \\ \hline
{\rule[0mm]{0mm}{4.5mm}\mbox{{\bf m)}}} &  \widetilde{A} \cap N = \Omega_9(q).a, \widetilde{B} \cap N = N_1, a \le 2, m = 8 & s & r \\ \hline
{\rule[0mm]{0mm}{4.5mm}\mbox{{\bf n)}}} & N = \Omega_{24}^+(2), \widetilde{A} \cap N = Co_1, \widetilde{B} \cap N = N_1 & 13 & 683 \\ \hline
\end{array}
\]\vspace{-.4cm}
\caption{Table 13}\label{tog+}
\end{table}

We treat now the case $m = 4$ and  consider the maximal factorizations according to \cite{LPS}, Tables 4. Here:

\begin{itemize}\item[-] The primes $r$ and $s$ ($r = 7$ if $q = 2$) are independent.\end{itemize}
For all maximal factorizations $G = \widetilde{A}\widetilde{B}$
with $\Omega_8^+(q) \le G \le \mbox{Aut}(\Omega_8^+(q))$ listed in Table 4 of \cite{LPS}, $r$ divides $|A \cap N|$ and
$s$ divides $|B \cap N|$, or vice versa.

\subsection{Exceptional groups of Lie type}

For the exceptional groups of Lie type all factorizations (not only the maximal ones) have been determined in \cite{HLS};
see also Table 5 in \cite{LPS}.

\begin{itemize}
\item[{\bf a)}] $N = G_2(q)$, $q = 3^e$, $A \cap N = SL_3(q)$ or $SL_3(q).2$, $B \cap N = SU_3(q)$ or $SU_3(q).2$:\\
By Table 3 of \cite{AI} (with corrected misprint), there exist independent prime divisors $r$ of
$q^2+q+1 \, | \, |A \cap N|$ and $s$ of $q^2-q+1 \, | \, |B \cap N|$.

\item[{\bf b)}] $N = G_2(q)$, $q = 3^e$, $e$ odd, $A \cap N = SL_3(q)$ or $SL_3(q).2$, $B \cap N = \, ^2G_2(q)$:\\
This case is handled as in {\bf a)}.

\item[{\bf c)}] $N = G_2(4)$, $A \cap N = J_2$, $B \cap N = SU_3(4)$ or $SU_3(4).2$:\\
7 divides $|A \cap N|$ and 13 divides $|B \cap N|$. 7 and 13 are independent.

\item[{\bf d)}] $N = G_2(4)$, $A \cap N = (G_2(2) \times 2)\cap N$, $B \cap N = (SU_3(4).4) \cap N$:\\
This case is handled as in {\bf c)}.

\item[{\bf e)}] $N = F_4(q)$, $q = 2^e$, $A \cap N = Sp_8(q)$, $B \cap N = {}^3D_4(q)$ or $^3D_4(q).3$:\\
By Table 3 of \cite{AI}, there exist independent prime divisors $r$ of
$q^4+1 \, | \, |A \cap N|$ and $s$ of $q^4-q^2+1 \, | \, |B \cap N|$.

\end{itemize}

\subsection{Sporadic groups}\label{spor}

We use the list of maximal factorizations in Table 6 of \cite{LPS}. \begin{itemize}\item[-] In all cases there are
independent primes dividing $|A \cap N|$ and $|B \cap N|$, respectively. They are presented Table~\ref{tsg}.\end{itemize}
Note that Out$(N) = 1$ for $N \cong M_{11}, M_{23}, M_{24}, Ru, Co_1$.
\begin{table}[H]
\centering
\[ \begin{array}{|l|l|c|c|} \hline
{\rule[0mm]{0mm}{4.5mm}\mbox{{\bf a)}}} & N = M_{11}, \widetilde{A} = L_2(11), \widetilde{B} = M_{10} \mbox{ or } M_9.2 & 11 & 3 \\ \hline
{\rule[0mm]{0mm}{4.5mm}\mbox{{\bf b)}}} & N = M_{12}, \widetilde{A} \cap N = M_{11}, \widetilde{B} \cap N = M_{11}, L_2(11), 11 \not\:\mid |A \cap N| & 3 & 11 \\ \hline
{\rule[0mm]{0mm}{4.5mm}\mbox{{\bf c)}}} & N = M_{12}, \widetilde{A} \cap N = M_{11}, \widetilde{B} \cap N = M_{11}, L_2(11), 11 \, | \, |A \cap N| & 11 &3 \\ \hline
{\rule[0mm]{0mm}{4.5mm}\mbox{{\bf d)}}} & N = M_{12}, \widetilde{A} \cap N = M_{11}, \widetilde{B} \cap N = M_{10}.2, M_9.S_3, 2 \times S_5, & 11 & 3 \\
 & 4^2 \times D_{12} \mbox{ or } A_4 \times S_3 &  &  \\ \hline
{\rule[0mm]{0mm}{4.5mm}\mbox{{\bf e)}}} & N = M_{12}, \widetilde{A} \cap N = L_2(11), \widetilde{B} \cap N = M_{10}.2 \mbox{ or } M_9.S_3 & 11 & 3 \\ \hline
{\rule[0mm]{0mm}{4.5mm}\mbox{{\bf f)}}} & N = M_{22}, \widetilde{A} \cap N = L_2(11).2 \cap N, \widetilde{B} \cap N = L_3(4).2 \cap N & 11 & 7 \\ \hline
{\rule[0mm]{0mm}{4.5mm}\mbox{{\bf g)}}} & N = M_{23}, \widetilde{A} = 23.11,  \widetilde{B} = M_{22}, M_{21}.2 \mbox{ or } 2^4.A_7 & 23 & 7 \\ \hline
{\rule[0mm]{0mm}{4.5mm}\mbox{{\bf h)}}} & N = M_{24}, \widetilde{A} = M_{23}, \widetilde{B} = L_2(23), 23 \, | \, |A| & 23 & 2 \\ \hline
{\rule[0mm]{0mm}{4.5mm}\mbox{{\bf i)}}} & N = M_{24}, \widetilde{A} = M_{23}, \widetilde{B} = L_2(23), 23 \not\:\mid |A| & 2 & 23 \\ \hline
{\rule[0mm]{0mm}{4.5mm}\mbox{{\bf j)}}} & N = M_{24}, \widetilde{A} = M_{23}, \widetilde{B} = M_{12}.2, 2^6.3.S_6, L_2(7) \mbox{ or } 2^6(L_3(2)\times S_3) & 23 & 2 \\ \hline
{\rule[0mm]{0mm}{4.5mm}\mbox{{\bf k)}}} & N = M_{24}, \widetilde{A} = L_2(23), \widetilde{B} = M_{22}.2, 2^4.A_8 \mbox{ or } L_3(4).S_3 & 23 & 2 \\ \hline
{\rule[0mm]{0mm}{4.5mm}\mbox{{\bf l)}}} & N = J_2, \widetilde{A} \cap N = U_3(3), \widetilde{B} \cap N = A_5 \times D_{10} & 7 & 5 \\ \hline
{\rule[0mm]{0mm}{4.5mm}\mbox{{\bf m)}}} & N = J_2, \widetilde{A} \cap N = U_3(3).2 \cap N, \widetilde{B} \cap N = 5^2(4 \times S_3) \cap N & 7 & 5 \\ \hline
{\rule[0mm]{0mm}{4.5mm}\mbox{{\bf n)}}} & N = HS, \widetilde{A} \cap N = M_{22}, \widetilde{B} \cap N = U_3(5).2, 7 \, | \, |A \cap N| & 7 & 5 \\ \hline
{\rule[0mm]{0mm}{4.5mm}\mbox{{\bf o)}}} & N = HS, \widetilde{A} \cap N = M_{22}, \widetilde{B} \cap N = U_3(5).2, 7 \not\:\mid |A \cap N| & 11 & 7 \\ \hline
{\rule[0mm]{0mm}{4.5mm}\mbox{{\bf p)}}} & N = HS, \widetilde{A} \cap N = M_{22}.2 \cap N, \widetilde{B} \cap N = 5^{1+2}.2^5 \cap N & 7 & 5 \\ \hline
\end{array}
\]
\end{table}
\begin{center}\vspace{-.5cm} Table 14 (contd.)\end{center}

\begin{table}[H]
\centering
\[ \begin{array}{|l|l|c|c|} \hline
{\rule[0mm]{0mm}{4.5mm}\mbox{{\bf q)}}} & N = He, \widetilde{A} \cap N = Sp_4(4).2, \widetilde{B} \cap N = 7^2.SL_2(7) & 17 & 7 \\ \hline
{\rule[0mm]{0mm}{4.5mm}\mbox{{\bf r)}}} & N = He, \widetilde{A} \cap N = Sp_4(4).4 \cap N, \widetilde{B} \cap N = 7^{1+2}(S_3 \times 6) \cap N & 17 & 7 \\ \hline
{\rule[0mm]{0mm}{4.5mm}\mbox{{\bf s)}}} & N = Ru, \widetilde{A} = L_2(29), \widetilde{B} = {}^2F_4(2) & 29 & 13 \\ \hline
{\rule[0mm]{0mm}{4.5mm}\mbox{{\bf t)}}} & N = Suz, \widetilde{A} \cap N = G_2(4) \cap N, \widetilde{B} \cap N = U_5(2) \mbox{ or } 3^5.M_{11} & 13 & 11 \\ \hline
{\rule[0mm]{0mm}{4.5mm}\mbox{{\bf u)}}} & N = Fi_{22}, \widetilde{A} \cap N = {}^2F_4(2)', \widetilde{B} \cap N = 2.U_6(2) & 13 & 11 \\ \hline
{\rule[0mm]{0mm}{4.5mm}\mbox{{\bf v)}}} & N = Co_1, \widetilde{A} = Co_2 \mbox{ or } Co_3, \widetilde{B} \cap N = 3.Suz.2 \mbox{ or } (A_4 \times G_2(4)).2 & 23 & 13 \\ \hline
\end{array}
\]\vspace{-.4cm}
\caption{Table 14}\label{tsg}
\end{table}

\end{itemize}

This completes the proof of Theorem \ref{as}.\bigskip

As a consequence of the proof we state the following result:

\begin{te}\label{fact}
Let $N$ be a simple group of Lie type (classical or exceptional) or a sporadic simple group.\\
Let $G = AB$ be a factorized almost simple group with socle $N$, i.e. $N \le G \le \mbox{\upshape Aut}(N)$, $A \le G$,
$B \le G$ with $N \not\le A$, $N \not\le B$.
\begin{itemize}
\item[a)] If $N \not\cong PSp_6(2)$ and $U_4(3)$, there exist independent primes with respect to $N$, one
dividing $|A \cap N|$, the other one dividing $|B \cap N|$.\smallskip

If $N \cong PSp_6(2)\cong \Omega_7(2)$, there exists a (non-maximal) factorization $N = AB$
of type g), Section \ref{symp}, such that there are no independent primes dividing $|A|$ and $|B|$, respectively.\smallskip

If $N \cong U_4(3)\cong P\Omega_6^{-} (2)$, then for each of the three factorization-types in m), Section \ref{un}, there exists a
(non-maximal) factorization $G = AB$ with $|G:N| = 2$ such that there are no independent primes dividing
$|A \cap N|$ and $|B \cap N|$, respectively.
\item[b)] If in addition to the two exceptional cases in a) $N \not\cong L_2(8),\ L_2(p),$ $p$ Mersenne prime, $L_3(3),\ L_4(2)$,
and $U_3(3)$, then there exist primes dividing $|A \cap N|$ and $|B \cap N|$ that are independent with
respect to $\mbox{\upshape Aut}(N)$.\smallskip

For $N \cong L_2(8),\ L_2(p),$ $p$ Mersenne prime, $L_3(3)$, or $U_3(3)$, all pairs of primes dividing $|N|$ are
not independent with respect to $\mbox{\upshape Aut}(N)$.\smallskip

 For $N \cong L_4(2)$ there exist factorizations $G = AB$ without primes
independent with respect to $\mbox{\upshape Aut}(N)$ dividing $|A \cap N|$ and $|B \cap N|$, respectively.
\end{itemize}
\end{te}

{\bf Proof.} 
 \ {\it a)} We mention that we can replace $G$ by $AN \cap BN = (A \cap BN)(B \cap AN)$
to satisfy the condition $G = AN = BN$ in Theorem \ref{as}; since $N \not\le A$ and
$N \not\le B$, none of the factors $(A \cap BN)$ and $(B \cap AN)$ equals $AN \cap BN$.\smallskip

The first assertion follows now from an inspection
of the proof of Theorem \ref{as}. Note however that in \cite{LPS} the group $N = L_4(2) \cong A_8$ is treated
as alternating group. Here the primes 2 and 7 are independent with respect to $N$, and any pair of a 2-Sylow
subgroup with a 7-Sylow subgroup generates $N$. This proves the assertion for $L_4(2)$.\smallskip

By Table 1 of \cite{LPS} there is a factorization $PSp_6(2) = \widetilde{A}\widetilde{B}$ with
$\widetilde{A} = O_6^-(2) = U_4(2):2$, $\widetilde{B} = O_6^+(2)$ and $\widetilde{A} \cap \widetilde{B} = S_6 \times 2$ is
maximal in $\widetilde{A}$ and $\widetilde{B}$. We infer from Table 1 of \cite{LPS} again that there is a
factorization $\widetilde{A} = A(\widetilde{A} \cap \widetilde{B})$ with $A = GU_3(2)$ of order $2^3 \cdot 3^4$. It
follows that $PSp_6(2) = A\widetilde{B}$. Since an element of order 5 in $PSp_6(2)$ is normalized by $S_3.4$ and
an element of order 7 by a cyclic group of order 6, none of the prime divisors of $|A|$ is independent of the
prime divisors of $|\widetilde{B}|$.

Therefore the exceptional case $PSp_6(2)$ does actually occur.\smallskip

Let $G = N\langle x\rangle$ with $N = U_4(3)$ and $x$ the square of the diagonal automorphism.
There are two conjugacy classes
of involutions, represented by 2B and 2C in \cite{Atl}, such that $G =N\langle s\rangle = N\langle t\rangle$,
$s$ of type 2B, $t$ of type 2C. There exists a subgroup $A_0 \cong L_3(4)$ in $N$ normalized by $s$ and a 3-Sylow subgroup
$D$ of $N$ normalized by $t$ such that $st$ generates together with a 2-Sylow subgroup of $A_0$ a 2-Sylow subgroup
of $N$. Setting $A = A_0\langle s \rangle$ and $B = D\langle t \rangle$, then $G = AB$, $A \cap N = A_0$, $B\cap N = D$. The prime 3 is not independent
of the other prime divisors of $|N|$: the normalizer of an element of order 7 has order 21 and there is a subgroup
of order $3^4$ normalized by $A_6$, so in particular by elements of order 5 and 2.

Therefore the exceptional case $N = U_4(3)$ with $|G:N| \ge 2$ does
actually occur.\medskip

{\it b)} For $N = L_4(2)$ a factorization $N = AB$ with $A \cong A_7$ and $B$ a 2-Sylow subgroup provides
an example: 2 is not independent of any other prime with respect to Aut$(N)$.\smallskip

For all other groups, apart from those already treated as exceptions in a), only such cases in the proof of
Theorem \ref{as} have to be considered where one of the independent primes divides $|\mbox{\upshape Out}(N)|$.
This leads to the groups $N \cong L_2(8), L_2(p),$ $p$ Mersenne prime, $L_3(3), U_3(3)$. It is easily checked
that all pairs of primes dividing $|N|$ are not independent with respect to Aut$(N)$.\qed
\bigskip

\begin{rem}
We have excluded the alternating groups from consideration in Theorem \ref{fact} and add some remarks about them here. \smallskip\\
We recall that $A_5\cong L_2(5)$, $A_6\cong L_2(9)$ and $A_8\cong L_4(2)$.\smallskip\\
More generally, there are many examples of alternating groups $A_n$ that possess a factorization $A_n = AB$
without independent primes dividing $|A|$ and $|B|$, respectively. For instance, if $n = 2^m$, $2^m - 1$ not prime,
or $n = 3^m$, $m \ge 2$, let $A = A_{n-1}$ and $B$ a 2-Sylow subgroup of $A_n$ in the first case and
a 3-Sylow subgroup in the second case. Then for $n = 2^m$ and every prime $p$, $2 < p < n$, a $p$-cycle
is normalized by a suitable involution in $A_n$ since $p < n-2$. If $n = 3^m$ and if $n-2$ is not a prime, then for
every prime $p$, $p \ne 3$, a $p$-cycle for $p$ odd or a product of two disjoint 2-cycles is
centralized by a suitable 3-cycle. If $q:= n-2$ is a prime, then 3 divides $q-1$ and a $q$-cycle
is normalized by a suitable element of order 3.\smallskip\\
On the other hand, if $p \ge 5$ is a prime, then every non-trivial factorization $A_p = AB$
yields independent primes dividing the order of the factors. This can be seen as follows: In $A_5$ the primes 3 and 5
are independent and a 3-Sylow together with a
5-Sylow subgroup generate the whole group; the same is true in $A_7$ with the primes 5 and 7. This proves the
assertion for $p = 5, 7$. For $p \ge 11$ there exists a prime $q$ with $\frac{p}{2} \le q \le p-4$; one can choose
$q = 7$ for $p = 11 $ or 13 and for $p \ge 17$ this follows from \cite{Ra}. By \cite{LPS}, Theorem D, we may assume that
a $p$-cycle is contained in $B$ and that $A_{p-4} \le A$ since there exists no 5-transitive group of prime degree $p$
other than $A_p$ or $S_p$. Then a $q$-cycle is contained in $B$.
As $q$ does not divide $p-1$, $p$ and $q$ are independent.\\
There are also other examples for both situations. But it appears to be difficult to determine all alternating groups where
every non-trivial factorization yields independent primes dividing the order of the factors.
\end{rem}

\section{The general case}

{\bf Proof of Main Theorem}\\
That (1) implies (2) is trivial. Also that (3) implies (1) is obvious: If $G = AB$ with $[A,B] \le G_{\So}$,
then $\langle a,b \rangle' \le G_{\So}$ whence $\langle a,b \rangle$ is soluble for all $a \in A$ and $b \in B$.\smallskip

It remains to prove that (2) implies (3). Suppose not and let $G$ be a minimal counterexample.
The proof proceeds in several steps.\smallskip

(i) $G$ has a unique minimal normal subgroup $N$, $N$ is non-abelian, $C_G(N) = 1$ (in particular, $G_{\So} = 1$),
$N \leq [A,B] \leq R$ where $R/N = (G/N)_{\So}$:\smallskip

If $N$ is a minimal normal subgroup of $G$, then one verifies easily that also $G/N = AN/N \cdot BN/N$
satisfies (2).
$G/N$ is not a counterexample whence $[A,B] \leq R$ with $R/N = (G/N)_{\So}$. Since
$[A,B] \trianglelefteq G$ and $[A,B] \neq 1$ ($G$ being a counterexample), $N \leq [A,B]$.
Clearly $N$ is non-abelian, for otherwise $[A,B]$ would be soluble.\\
Suppose there are two distinct minimal normal subgroups $N_1, N_2$ of $G$. Let $R_i/N_i = (G/N_i)_{\So}$,
$i = 1,2$. Choose $d$ such that $R_i^{(d)} = N_i$, $i = 1,2$ ($d$-th derived subgroup). Then
$[A,B]^{(d)} \leq R_1^{(d)} \cap R_2^{(d)} = N_1 \cap N_2 = 1$ and $[A,B]$ is soluble, contradiction.\\
Therefore $G$ has a unique minimal normal subgroup $N$. As $N$ is non-abelian,
$N \cap C_G(N) = 1$, and it follows from $C_G(N) \trianglelefteq G$ that $C_G(N) = 1$.\smallskip

(ii) $N \not\le A$ and $N \not\le B$:\smallskip

Suppose that $N \le A$. Since $B \ne 1$, we can choose $b \in B$ of prime order, say $p$. Then $N\langle b \rangle =
(A \cap N\langle b \rangle)\langle b \rangle$ satisfies (2).
If $N\langle b \rangle < G$, then $[N,\langle b \rangle] \le (N\langle b \rangle)_{\So} \le C_{N\langle b \rangle}(N) = 1$
by (i). Therefore $b \le C_G(N) = 1$, contradiction. Hence $N\langle b \rangle = G$.\\
Suppose first that $N$ is simple. If $G = A$, then Theorems 1 and 2 of \cite{GL} imply that $b \in G_{\So} = 1$,
contradiction. If $G = A\langle b \rangle > A = N$, let $1 \neq g \in G$ be a $q$-element for a prime $q \ne p$. Then
$g \in N = A$ and $\langle b,g \rangle$ is soluble by (2). Again the same theorems as before imply that
$b \in G_{\So} = 1$, a contradiction.\\
Hence $N$ is not simple and in particular $b \not\in N$. Consequently, $N$ is the direct product
of $p$ copies of a non-abelian simple group $L$, permuted transitively by $b$.
Therefore we can assume w.l.o.g. that there are automorphisms $\alpha_i$ of $L$, $1,\ldots, p$, such that $(l_1,\ldots,l_p)^b =
(\alpha_2(l_2),\ldots,\alpha_p(l_p),\alpha_1(l_1))$ for all $(l_1,\ldots,l_p) \in N$.\\
Assume first that $p = 2$. Let $q \ge 5$ be a prime divisor of $|L|$. Since $G_{\So} = 1$, by \cite{GGKP1}, \cite{GGKP2} or
\cite{Gu} there exist (conjugate) elements $m_1, m_2 \in L$ of order $q$ such that $\langle m_1,m_2 \rangle$ is not soluble.
Set $n = (m_1,\alpha_2^{-1}(m_2))$.
Then $n$ is of order $q$ and $\langle n,b \rangle \ge \langle n, n^b \rangle = \langle (m_1,\alpha_2^{-1}(m_2)), (m_2, \alpha_1(m_1)\rangle$
is not soluble, contradicting (2).\\
Now let $p \ge 3$. By Theorem A of \cite{MSW} there exist three involutions
$m_1, m_2, m_3$ in $L$ such that $\langle m_1, m_2, m_3\rangle$ is not soluble unless $L = U_3(3)$. In the latter
case we let $m_1, m_2$ be elements of order 8 that generate a maximal subgroup $3^{1+2}:8$ of $L$ and choose $m_3$
as a 2-element outside this maximal subgroup. Setting
$l_1 = m_1$, $l_2 = \alpha_2^{-1}(m_2)$, $l_3 = \alpha_3^{-1}\alpha_2^{-1}(m_3)$ and $l_i = 1$ for $3 < i \le p$ in case $p \ge 5$,
the element $n = (l_1,\ldots, l_p) \in N$ is
a 2-element. Then conjugating $n$ with $b$, $b^2$, it follows as above that $\langle n,b \rangle$ is not
soluble, contradicting (2).\\
Hence $N \not\le A$ and analogously $N \not\le B$.\smallskip

(iii) $AN = BN = G$ and $G/N$ is soluble.\smallskip

Suppose $AN = A(AN \cap B) < G$. Then $[A, AN \cap B] \le (AN)_{\So} \le C_G(N) = 1$. It follows that $AN \cap B$
is normal in $AN$ and that $N \cap B$ is normalized
by $B$ and centralized by $A$, whence $N \cap B = 1$ or $N \cap B = N$. The latter is impossible by (ii).
Then $[N, AN \cap B] \le N \cap AN \cap B = N \cap B = 1$ whence $AN \cap B \le C_G(N) = 1$. This implies
$AN = A$, contradicting (ii).\\
This shows that $AN = G$ and analogously $BN = G$.\\
Because of $AN = BN = AB = G$ and the fact that $A$ and $B$ satisfy (2) it follows that all for all primes $p \ne q$
every $p$-element and every $q$-element of $G/N$ generate a soluble group. Then by Theorem B of \cite{DGHP} (applied
to a hypothetical non-cyclic composition factor of $G/N$) or Theorems 1 and 2 of \cite{GL} it follows that $G/N$ is
soluble.\smallskip

(iv) Let $N$ be the direct product
of $k$ copies of a non-abelian simple group $L$. We denote the $i$-th component by $L_i$.
Since $G$ is a counterexample, it follows from
(i), (ii), (iii) and Theorem \ref{as} that $k \ge 2$.\medskip

 Let $M_i =$ Aut$(L_i) \cong$ Aut$(L)$ and identify $L_i$ with Inn$(L_i)$ whence
$L_i \le M_i$, $i = 1,\ldots,k$. Then $N = L_1 \times \cdots \times L_k \le G \le (M_1 \times \cdots \times M_k)T
\cong$ Aut$(L)\, wr\, S_k$ where $T \cong S_k$, the symmetric group of degree $k$. Set $M = M_1
\times \cdots \times M_k$. We denote the elements of $T$ by permutations $t$ where $t$ acts on $M$
via $(m_1,\ldots,m_k)^t = (m_{1t^{-1}},\ldots,m_{kt^{-1}})$. \\
Let $\rho_i$ denote the projection of $M$ onto $M_i$, $1\le i \le k$.\\
For $X = G,A,B$ and $1 \le i \le k$ set\\ $W_{i,X} = \{(m_1,\ldots,m_k)t \in X \,|\, m_j \in M_j \mbox{ and } t \in T \mbox{ fixes } i\}$
and\\
$V_{i,X} = \{m_i \in M_i \,|\, \mbox{there exists }(m_1,\ldots,m_k)t \in W_{i,X}\}$.\\
It is clear that $W_{i,X}$ is a subgroup of $X$ and $X \cap M \le W_{i,X}$. Also $V_{i,X}$ is a subgroup of $M_i$,
$\rho_i(X \cap M) \le V_{i,X}$.\smallskip

We claim:\medskip

(v) For each $1 \le i \le k$, $V_{i,G} = V_{i,A}$ and $V_{i,B} = 1$ or vice versa.\smallskip

It suffices to prove the assertion for $i = 1$ (the other cases being analogous) and we set $V_X$ and $W_X$ for
$V_{1,X}$ and $W_{1,X}$, respectively.\\
Let $l_1 \in L_1$ be arbitrary. There exist $a \in A$ and $b \in B$ with $(l_1,1,\ldots,1) = ab$. \\
Let $a = (m_1,\ldots,m_k)t$ and assume w.l.o.g. that in the decomposition of $t$ in disjoint cycles the cycle
containing 1 is just $(1,\ldots,d)$ for some $1 \le d \le k$. Clearly, $b = t^{-1}(m_1^{-1}l_1,m_2^{-1},\ldots,m_k^{-1})$.\\
Set $m = (m_1,\ldots,m_k)$ and $m' = (m_1^{-1}l_1,m_2^{-1},\ldots,m_k^{-1})$.\\ Then
$a^d = m\cdot m^{t^{-1}}\cdots m^{t^{-(d-1)}}t^d = (m_1m_2\cdots m_d,*,\ldots,*)t^d\in W_A$ and\\
$b^d = m'^t \cdots m'^{t^d}t^{-d} = (m_d^{-1}m_{d-1}^{-1}\cdots m_2^{-1}m_1^{-1}l_1,*,\ldots,*)t^{-d} \in W_B$.\\
It follows that $(l_1,*,\ldots,*) = a^db^d \in W_AW_B$.
Therefore, $l_1 \in V_AV_B$, whence $L_1 \subseteq V_AV_B$.\\
Let $g \in W_G$ be arbitrary. By (2) $g = an$, $a = (m_1,\ldots,m_k)t \in W_A$, $n = (l_1,\ldots,l_k)\in N$. By what we have proved before,
$m_1l_1 \in V_AV_B$, whence $V_G = V_AV_B$.\\
We claim that $V_A$ and $V_B$ satisfy (2) in $V_G$. Since $W_A$ and $W_B$ satisfy (2), it suffices to show that for
every prime $p$ and every $p$-element $m \in V_X$ there exists a $p$-element $(m',m_2,\ldots,m_k)t \in W_X$ ($X = A, B$)
with $\langle m \rangle = \langle m' \rangle$. This follows easily: Choose $y = (m,m'_2,\ldots, m'_k)\tilde{t} \in W_X$. Let
the order of $y$ be $p^a\cdot r$ where $p$ does not divide $r$. Then $y^r = (m^r,\ldots){\tilde{t}}^r$
since $t$ fixes position 1. The assertion follows with $m' = m^r$.\\
Minimality of $G$ (recall that $k \geq 2$) now yields that $[V_A,V_B] \le (V_G)_{\So}$. \\
Since $L_1 \le V_G \le M_1$ and since $L_1$ is the unique minimal normal subgroup of $M_1$, it follows that
$[V_A,V_B] = 1$. This means that one of $V_A$ and $V_B$ contains $L_1$ and the other is trivial. Hence (v) holds.\smallskip

(vi) $G \cap M = N$; moreover with appropriate choice of notation $B \cap M = 1$ and $\rho_i(A\cap N) = L_i$ for $1 \le i \le k$.\smallskip

By (v) we may assume that $V_{1,G} \times \cdots \times V_{k,G} = V_{1,A} \times \cdots \times V_{h,A} \times
V_{h+1,B} \times \cdots \times V_{k,B}$ and $V_{1,B} = \cdots V_{h,B} = V_{h+1,A} = \cdots V_{k,A} = 1$
for some $0 \le h \le k$. W.l.o.g. assume that $h \ge 1$. Then $L_1 \le V_{1,G} = V_{1,A}$.\\
By (iii) $G/N$ is soluble. Then
$A/(A\cap N)$ is soluble and taking iterated derived groups, it follows
that $L_1 \le \rho_1(A \cap N)$. Since $G = AN$ by (iii),
$A$ acts transitively on $\{L_1,\ldots, L_k\}$. As $A\cap N$ is normal in $A$, $N \le \rho_1(A\cap N) \times
\cdots \times \rho_k(A\cap N)$. Moreover, $B\cap M = 1$. Since $BN = G$ by (iii), the latter implies $G\cap M = N$.\smallskip

(vii) Final contradiction.\smallskip

It is well known (see for instance \cite{BE}, Proposition 1.1.39) that
because of (vi) there is a partition of $\{1,\ldots,k\}$ into subsets $J_1,\ldots,J_d$ such that
$A \cap N = \vartimes_{i=1}^d A_i$ where $A_i = A \cap \vartimes_{j \in J_i} L_j$ and $|A_i| = |L|$, $1 \le i \le d$.\\
Suppose that $A \cap L_i = 1$ for all $1 \le i \le k$. Then $d \le k/2$ (note that $k \ge 2$, see (iv)),
whence $|A\cap N| \le |L|^{k/2}$. Since $G = AN$,
$|G| = |A||N|/|A\cap N| \ge |A||L|^{k/2}$. Since $B\cap N = 1$ by (5), $B$ is isomorphic to a soluble
subgroup of $S_k$ and hence $|B| \le 3^{k-1}$ by a result of Dixon (see \cite{DM},
Theorem 5.8B). Therefore $|G| \le |A||B| \le |A|\cdot3^{k-1}$. It follows that $|L|^{k/2} \le 3^{k-1} < 9^{k/2}$,
whence $|L| < 9$, a contradiction.\\
Therefore we may assume that $A \cap L_1 \ne 1$. If $(l,1,\ldots,1) \in A \cap L_1$ with $l \ne 1$, then
because of $\rho_1(A\cap N) = L_1 (\cong L)$ all $([l,l'],1,\ldots,1), l' \in L$ are in $A \cap L_1$. But the elements
$[l,l']$ generate $L$, whence $A \cap L_1 = L_1$.\\
As $A$ acts transitively on $\{L_1,\ldots, L_k\}$, it follows that $N \le A$ which
contradicts (ii).\qed

\begin{cor}\label{cote1} Let $N$ be a non-abelian simple group and $N \le G = AB \le \mbox{\upshape Aut}(N)$ with
subgroups $A$ and $B$ satisfying (2) of the Main Theorem, then $A = 1$ or $B = 1$.
\end{cor}

{\bf Proof.} By the Main Theorem, $[A,B] \le G_{\So} = 1$. Therefore $A$ and $B$ are normal in $G$, whence
$A = 1$ or $B = 1$ or $N \le A \cap B$. But the latter implies $N = N' \le [A,B] = 1$, a contradiction.\qed \\

There is another corollary to the Main Theorem that generalizes a result of Carocca \cite{C} on the
solubiltity of $\So$-connected products of finite soluble groups:

\begin{cor}\label{corCa} Let $G = AB$ be a finite group with $\So$-connected subgroups $A, B$. Then $A_{\So} = A \cap G_{\So}$
and $B_{\So} = B \cap G_{\So}$.\\
In particular, if $A$ and $B$ are soluble then $G$ is soluble.
\end{cor}

{\bf Proof.} By the Main Theorem, $A_{\So}G_{\So}$ is a normal subgroup of $G$. Since $A_{\So}G_{\So}/G_{\So}$ is
soluble, $A_{\So}G_{\So}$ is soluble. Consequently, $A_{\So} \le G_{\So}$ whence $A_{\So} = A \cap G_{\So}$.\qed

\bigskip

P. HAUCK\\
{\small Fachbereich Informatik,
Universit\"{a}t T\"{u}bingen,\\
Sand 13, 72076 T\"{u}bingen, Germany\\
e-mail: peter.hauck@uni-tuebingen.de}\\
 \\
 L.~S. KAZARIN\\
{\small Department of Mathematics, Yaroslavl P. Demidov State University\\
Sovetskaya Str 14, 150014 Yaroslavl, Russia\\
e-mail: Kazarin@uniyar.ac.ru}\\
 \\
A. MART\'{I}NEZ-PASTOR\\
{\small Escuela T\'{e}cnica Superior de  Ingenier\'{\i}a Inform\'{a}tica,\\
Instituto Universitario de Matem\'{a}tica Pura y  Aplicada IUMPA\\
Universidad Polit\'{e}cnica de Valencia,
 Camino de Vera, s/n,  46022 Valencia, Spain\\
e-mail: anamarti@mat.upv.es}\\
 \\
M.~D. P\'{E}REZ-RAMOS\\
{\small Departament de Matem\`{a}tiques, Universitat de Val\`{e}ncia,\\
C/ Doctor Moliner 50, 46100 Burjassot
(Val\`{e}ncia), Spain\\
e-mail: Dolores.Perez@uv.es}

\end{document}